\documentclass[a4paper,11pt]{article}
\usepackage[T1]{fontenc}
\usepackage{lmodern,amsmath,amssymb,amsthm,mathtools,empheq}
\allowdisplaybreaks[4]
\usepackage[a4paper,left=25mm,right=25mm,top=26mm,bottom=25mm,headheight=14pt]{geometry}
\usepackage{fancyhdr}
\usepackage{xcolor}
\usepackage[bookmarksnumbered,colorlinks=true,linkcolor=blue!50!black,citecolor=blue!50!black,urlcolor=blue!50!black]{hyperref}
\hypersetup{pdftitle={Global strong solutions to the 2D compressible Navier--Stokes equations with vacuum and large initial data}}

\numberwithin{equation}{section}
\newtheorem{theorem}{Theorem}[section]
\newtheorem{lemma}[theorem]{Lemma}
\newtheorem{proposition}[theorem]{Proposition}
\theoremstyle{definition}

\newtheorem{remark}{Remark}[section]
\title{
Global strong solvability of the two-dimensional\\[2pt]
compressible Navier--Stokes equations for large initial data%
}

\author{Jie F{\small an}$^{a}$, Xiangdi H{\small uang}$^{b}$\thanks{Email addresses: fj0828@outlook.com (J. Fan); xdhuang@amss.ac.cn (X. D. Huang); muxilei@amss.ac.cn (M. X. Lei)}, Muxi L{\small ei}$^{b}$   \\  
	{\normalsize a. College of Science, Minzu University of China, Beijing 100081, P. R. China;}\\
    {\normalsize b. Institue of Mathematics, Academy of Mathematics and Systems Sciences,}\\
	{\normalsize Chinese Academy of Sciences, Beijing 100190, China;}\\
}
\date{}

\begin{document}

\begingroup
\makeatletter
\renewcommand{\@fnsymbol}[1]{}
\makeatother
\maketitle
\endgroup

\thispagestyle{plain}
\begin{abstract}
In 1995, Kazhikhov and Vaigant introduced a two-dimensional
barotropic compressible Navier--Stokes model with constant shear
viscosity and bulk viscosity $\lambda(\rho)=\rho^\beta$. For $\beta>3$,
they proved global existence of smooth solutions for arbitrarily
large initial data with strictly positive density. Subsequent work
allowed initial vacuum, widened the range of viscosity exponents,
and treated other domains [Jiu--Wang--Xin, J. Math. Fluid Mech. 16
(2014), 483--521; Huang--Li, J. Math. Pures Appl. 106 (2016), 123--154;
Huang--Li, SIAM J. Math. Anal. 54 (2022), 3192--3214; Fan--Li--Li,
Arch. Ration. Mech. Anal. 245 (2022), 239--278]. More recently, Fan,
-Jiang-Li-L\"u (arXiv:2611.20264, 2026) proved global strong
solvability for every $\beta>1$ in bounded simply connected domains
with Navier-slip boundary conditions. We prove global existence and uniqueness of strong solutions to
the Kazhikhov--Vaigant model on the torus $\mathbb{T}^2$ at
$\beta=1$, for every $\gamma>1$ and arbitrarily large initial data
allowing vacuum. To improve the $L^\gamma$ control furnished by
the basic energy estimate at the endpoint, we use the symmetry of
the periodic Riesz transforms to derive a new double-commutator
estimate. It yields a uniform $L^2$ bound for the density on each
finite time interval and a space--time integral bound for
$\rho^{\gamma+1}$. At $\beta=1$, the linear growth of the viscosity primitive is
insufficient to absorb the term with the factor $R\log R$ in the
effective viscous flux estimate, where $R(t)$ is the maximum density
up to time $t$. We retain the time-dependent weight $(R\log R)^{-1}$
inside the flux integral and reduce the density upper bound to
controlling the increase of $\log\log R$. The pressure contribution to the flux energy grows like
$\rho^{2\gamma-1}$ at large density. For $1<\gamma\le2$, the improved
density integrability makes the flux energy integrable in time.
For $\gamma>2$, we introduce a pressure correction that cancels
the mixed term linear in the flux and quadratic in the pressure
and yields pressure dissipation. On the time set where the pressure
contribution dominates the velocity dissipation, we use the positive
variation of the pressure integral to control the remaining growth
of $\log\log R$.
\end{abstract}
\noindent\textbf{MSC (2020)}\quad 35Q30, 35B45, 76N10.\par
\noindent\textbf{Keywords}\quad Compressible Navier--Stokes equations, global strong solutions,
density-dependent bulk viscosity, large initial data, vacuum.
\begingroup
\small
\tableofcontents
\endgroup

\section{Introduction and main results}\label{sec:1}
Let $\mathbb{T}^2=\mathbb{R}^2/\mathbb{Z}^2$. On $\mathbb{T}^2\times(0,\infty)$ we consider
the barotropic compressible Navier--Stokes equations
\begin{subequations}\label{eq:system}
\begin{empheq}[left=\empheqlbrace]{align}
&\rho_t+\operatorname{div}(\rho \mathbf{u})=0, \label{eq:continuity}\\
&(\rho \mathbf{u})_t+\operatorname{div}(\rho \mathbf{u}\otimes \mathbf{u})+\nabla P
 =\mu\Delta \mathbf{u}+\nabla\big((\mu+\lambda(\rho))\operatorname{div}\mathbf{u}\big).
 \label{eq:momentum}
\end{empheq}
\end{subequations}
Here $\rho$ is the density, $\mathbf{u}=(u_1,u_2)$ is the velocity, and $\mu>0$ is the constant
shear viscosity. We assume that the pressure and bulk viscosity satisfy
\begin{subequations}\label{eq:pressure-law}
\begin{align}
 P(\rho)      &= \rho^\gamma, \qquad \gamma>1, \label{eq:pressure-law-pressure}\\
 \lambda(\rho) &= \rho. \label{eq:pressure-law-viscosity}
\end{align}
\end{subequations}
We prescribe the initial density and momentum by
\begin{equation}\label{eq:initial-data}
 \rho(x,0)=\rho_0(x),\qquad (\rho \mathbf{u})(x,0)=m_0(x)=\rho_0(x)\mathbf{u}_0(x),
 \qquad x\in\mathbb{T}^2 .
\end{equation}

For the two-dimensional compressible Navier--Stokes equations with
constant shear viscosity and density-dependent bulk viscosity
$\lambda(\rho)=\rho^\beta$, Vaigant and
Kazhikhov~\cite{bib1} established the existence and uniqueness of global
strong solutions in rectangular domains when $\beta>3$ and the initial
density is bounded away from zero. Subsequently,
Jiu--Wang--Xin~\cite{bib2} established global classical solutions on $\mathbb{T}^2$
allowing initial vacuum under the same restriction on $\beta.$ Later, Huang and Li~\cite{bib3} developed commutator estimates
and lowered the threshold for global strong solvability to
$\beta>4/3$ for every $\gamma>1$, still allowing initial vacuum.
They also obtained a time-independent density upper bound and
convergence to equilibrium as $t\to\infty$ under the additional
conditions $\beta>3/2$ and $1<\gamma<4\beta-3$. For the Cauchy problem, Huang- Li~\cite{HuangLi2022} proved
the global existence and uniqueness of strong and classical solutions
for $\beta>4/3$ with a suitable decay rate of the initial density
(see also Jiu--Wang--Xin~\cite{JiuWangXin2013}). We can refer to Jiu--Wang--Xin~\cite{JiuWangXin2018} for positive
density at far field.
For the Navier-slip boundary condition,
Fan--Li--Li~\cite{FanLiLi2022}, using the Riemann mapping theorem and a
pull-back Green's function, established the existence and uniqueness
of strong solutions in bounded simply connected domains for
$\beta>4/3$. Later, Fan--Li--Wang~\cite{FanLiWang2023}
exploited the damping mechanism of the pressure to prove exponential
convergence to equilibrium. Their results also hold for periodic
domains. Fan-Jiang-Li-L\"u~\cite{bib7} extended global
strong solvability to $\beta>1$ under the boundary conditions
$u\cdot n=0$ and $\operatorname{curl}u=0$. Fan-Huang~\cite{bib4} established global strong solvability
with transported entropy for $\beta>4/3$ on periodic domains
and bounded domains with Navier-slip boundary conditions. For the radially symmetric free-boundary problem, Li-
Zhang~\cite{LiZhang2016} established the existence and uniqueness
of strong solutions for $\beta>1$ and estimated the rate of boundary
expansion when $1<\beta\le\gamma$. Huang--Meng--Ni~\cite{bib8}
used the geometric properties of two-dimensional radial symmetry
to extend the existence result to $\beta=1$. For the Dirichlet
problem, Huang-Su-Yan-Yu~\cite{HuangSuYanYu2024} provided
a representation of the effective viscous flux boundary term and
proved global strong solvability under radial symmetry
for $\beta>1$. They also obtained convergence to equilibrium when
$\max\{1,(\gamma+2)/4\}<\beta\le\gamma$. Huang and
Yan~\cite{HuangYan2024} removed the restriction $\beta\le\gamma$
when considering the asymptotic behavior of strong solutions
in the MHD system. Huang and Meng~\cite{bib9} studied the global
existence and asymptotic behavior of two-dimensional radially
symmetric strong solutions with swirl.

\subsection{Main theorem}\label{subsec:1.2}
We write
\[
 \begin{gathered}
 \int f\,dx=\int_{\mathbb{T}^2}f\,dx,\qquad
 \bar f=\frac{1}{|\mathbb{T}^2|}\int_{\mathbb{T}^2} f\,dx,\\
 L^r=L^r(\mathbb{T}^2),\qquad W^{s,r}=W^{s,r}(\mathbb{T}^2),\qquad H^s=W^{s,2}.
 \end{gathered}
\]
The material derivative is denoted by
\[
 \frac{D}{Dt}f=\dot f=f_t+\mathbf{u}\cdot\nabla f .
\]
\begin{theorem}\label{thm:1.1}
Assume that $\mu>0$ and that \eqref{eq:pressure-law} holds.
Suppose that, for some $q>2$, the periodic initial data $(\rho_0,m_0)$
satisfy
\begin{equation}\label{eq:initial-hypotheses}
 0\le \rho_0\in W^{1,q},\qquad
 \mathbf{u}_0\in H^1,\qquad
 m_0=\rho_0\mathbf{u}_0,\qquad
 M:=\int\rho_0\,dx>0.
\end{equation}
Then the problem \eqref{eq:system}--\eqref{eq:initial-data} has a unique
strong solution $(\rho,\mathbf{u})$ in $\mathbb{T}^2\times(0,\infty)$
satisfying
\begin{equation}\label{eq:strong-class}
\left\{
\begin{aligned}
 &\rho\in C([0,T];W^{1,q}),\qquad
   \rho_t\in L^\infty(0,T;L^2),\\
 &\mathbf{u}\in L^\infty(0,T;H^1)\cap L^2(0,T;H^2)
   \cap L^{(q+1)/q}(0,T;W^{2,q}),\\
 &t^{1/2}\mathbf{u}\in L^2(0,T;W^{2,q}),\qquad
   t^{1/2}\mathbf{u}_t\in L^2(0,T;H^1),\\
 &\rho^{1/2}\mathbf{u}_t\in L^2(0,T;L^2),\qquad
   \nabla\mathbf{u}\in L^1(0,T;L^\infty),
\end{aligned}
\right.
\end{equation}
for any $0<T<\infty$. Moreover, for each such $T$, there exists
a positive constant $C_T$, depending only on $T$, $\mu$, $\gamma$,
$\mathbb{T}^2$, $M$, $\|\rho_0\|_{L^\infty}$, and
$\|\mathbf{u}_0\|_{H^1}$, such that
\begin{equation}\label{eq:basic-regularity}
 \sup_{0\le t\le T}\|\rho(\cdot,t)\|_{L^\infty}\le C_T.
\end{equation}
\end{theorem}
\begin{remark}
The theory of global strong solutions for large initial data in the
Vaigant--Kazhikhov model has been developed in a series of works
\cite{bib1,bib2,bib3,bib6}. Theorem~\ref{thm:1.1} reaches the endpoint
$\beta=1$ on the torus for every $\gamma>1$, allowing initial vacuum.
The key step in the proof is to bound the density from above.
At the endpoint, the viscosity primitive \eqref{eq:viscosity-primitive}
has only linear growth, so an error of order $R\log R$ cannot be
absorbed directly. We use \eqref{eq:record-maximum} at the times when
$R$ increases. The weight in \eqref{eq:weighted-flux} allows us to
control the integral of $R'/(R\log R)$, which gives the density bound
\eqref{eq:basic-regularity} on each finite time interval.
\end{remark}
\subsection{Proof strategy}\label{subsec:1.3}
The key step in proving global existence is to establish an upper bound
for the density. At the linear bulk-viscosity endpoint, the primitive in
\eqref{eq:viscosity-primitive} grows linearly at large density. The estimate of Huang and Li~\cite{bib3} does not give
uniform control of density integrability as $\beta\downarrow1$. Let $R(t)$ be the maximum density up to time $t$ in \eqref{eq:density-maximum}, and let $\widetilde G$ be the
effective viscous flux in \eqref{eq:unshifted-flux}. At times when
$R$ increases, we have (see \eqref{eq:record-maximum})
\[
 \left(1+\frac{2\mu}{R}\right)R'+R^\gamma
 \le \|\widetilde G\|_{L^\infty}.
\]
On $[\tau,T_1]$, with a positive starting time $\tau$, we estimate
\begin{equation}\label{eq:intro-density-growth}
 \log\log R(T_1)-\log\log R(\tau)
 =\int_\tau^{T_1}\frac{R'(t)}{R(t)\log R(t)}\,dt.
\end{equation}
Proposition~\ref{prop:initial-interval} provides a positive starting
time $\tau$ at which $R(\tau)$, $\mathcal V(\tau)$, and
$\|\sqrt{\rho(\tau)}\,\dot{\mathbf u}(\tau)\|_{L^2}$ are bounded
in terms of the initial data, as in \eqref{eq:initial-interval-start}.
The monotonicity of $R$ determines the sign of the term arising from
differentiation of the weight. Choosing the interpolation exponents
in terms of $R(t)$, we obtain (see \eqref{eq:weighted-flux})
\[
 \int_\tau^{T_1}
 \frac{\|\widetilde G\|_{L^\infty}^{2}}
      {\mathcal V(R\log R)^2}\,dt
 \le C_{\tau,T}.
\]
We use \eqref{eq:record-maximum} to bound the increase of
$\log\log R$ in \eqref{eq:intro-density-growth}. Since $R(\tau)$
is bounded by \eqref{eq:initial-interval-start}, a bound on this
increase gives an upper bound for $R(T_1)$.

For $1<\gamma\le2$, the double-commutator argument yields the
space--time integral bound for $\rho^{\gamma+1}$ in
\eqref{eq:quadratic-density}. Consequently, the pressure contribution
to the flux energy is integrable in time. The Cauchy--Schwarz inequality
with \eqref{eq:subquadratic-energy-integral} and \eqref{eq:weighted-flux}
then bounds the integral in \eqref{eq:intro-density-growth}.

For $\gamma>2$, the pressure correction in \eqref{eq:corrected-energy}
cancels the term linear in the effective viscous flux and quadratic in
the pressure. The correction also produces pressure dissipation.
We split the time interval into two sets according to the relative
sizes of the pressure integral and the velocity dissipation. On the set where the
velocity dissipation dominates, the basic energy bound and the weighted
flux estimate control the increase of $\log\log R$. On the complementary
set, we use the pressure term in \eqref{eq:record-maximum} and the
positive-variation estimate \eqref{eq:pressure-positive-variation}.
The resulting density bound and the higher-order estimates in
Section~\ref{sec:4} support the continuation and vacuum approximation
in Section~\ref{sec:5}.

Section~\ref{sec:2} collects the preliminary results used below.
In Section~\ref{sec:3}, we establish a finite-time upper bound for
the density.  Section~\ref{sec:4} uses this bound to derive the
higher-order estimates. In Section~\ref{sec:5}, we prove
Theorem~\ref{thm:1.1} by continuing the positive-density solutions
and passing to the limit in the approximation.
\section{Preliminaries}\label{sec:2}
The following local existence result can be proved by adapting the
arguments in \cite{SalviStraskraba1993,Solonnikov1980}.
\begin{lemma}\label{lem:local-existence}
Assume that $(\rho_0,m_0)$ satisfies
\begin{equation}\label{eq:local-existence-data}
 \rho_0\in H^2,\qquad \mathbf{u}_0\in H^2,\qquad
 \inf_{x\in\mathbb{T}^2}\rho_0(x)>0,\qquad
 m_0=\rho_0\mathbf{u}_0.
\end{equation}
Then there are a small time $T_0>0$ and a constant $C_0>0$,
both depending only on $\|\rho_0\|_{H^2}$, $\|\mathbf{u}_0\|_{H^2}$, and
$\inf_{x\in\mathbb{T}^2}\rho_0(x)$ for fixed $\mu$, $\gamma$, and
$\mathbb{T}^2$, such that there exists a unique strong solution
$(\rho,\mathbf{u})$ to the problem
\eqref{eq:system}--\eqref{eq:initial-data} in
$\mathbb{T}^2\times(0,T_0)$ satisfying
\begin{equation}\label{eq:local-regularity}
\left\{
\begin{aligned}
 &\rho\in C([0,T_0];H^2),\qquad
   \rho_t\in C([0,T_0];H^1),\\
 &\mathbf{u}\in L^2(0,T_0;H^3),\qquad
   \mathbf{u}_t\in L^2(0,T_0;H^1),
\end{aligned}
\right.
\end{equation}
and
\begin{equation}\label{eq:local-density-lower-bound}
 \inf_{(x,t)\in\mathbb{T}^2\times(0,T_0)}\rho(x,t)
 \ge C_0>0.
\end{equation}
\end{lemma}

We shall use the following Poincar\'e--Sobolev and Brezis--Wainger
inequalities~\cite{LadyzenskajaSolonnikovUralceva1968,BrezisWainger1980,Engler1989}.
\begin{lemma}\label{lem:periodic-interpolation}
There exists a positive constant $C$ depending only on $\mathbb{T}^2$ such that every function
$u\in H^1(\mathbb{T}^2)$ satisfies for $2<p<\infty$,
\begin{equation}\label{eq:periodic-interpolation}
\begin{aligned}
 \|u-\bar{u}\|_{L^p}
 &\le Cp^{1/2}\|u-\bar{u}\|_{L^2}^{2/p}\|\nabla u\|_{L^2}^{1-2/p},\\
 \|u\|_{L^p}
 &\le Cp^{1/2}\|u\|_{L^2}^{2/p}\|u\|_{H^1}^{1-2/p}.
\end{aligned}
\end{equation}
Moreover, for $q>2$, there exists some positive constant $C$ depending only on $q$ and
$\mathbb{T}^2$ such that every function $v\in W^{1,q}(\mathbb{T}^2)$ satisfies
\begin{equation}
 \|v\|_{L^\infty}
 \le C\|\nabla v\|_{L^2}\ln^{1/2}\!\bigl(e+\|\nabla v\|_{L^q}\bigr)
     +C\|v\|_{L^2}+C.
\end{equation}
\end{lemma}

The following Poincar\'e type inequality can be found in
\cite{Feireisl2004}.
\begin{lemma}\label{lem:weighted-poincare}
Let $v\in H^1(\mathbb{T}^2)$, and let $\rho$ be a non-negative function such that
\[
 0<M_1\le\int_{\mathbb{T}^2}\rho\,dx,\qquad
 \int_{\mathbb{T}^2}\rho^\gamma\,dx\le M_2,
\]
with $\gamma>1$. Then there is a constant $C$ depending solely on $M_1$, $M_2$, and
$\gamma$ such that
\begin{equation}\label{eq:weighted-poincare-l2}
 \|v\|_{L^2(\mathbb{T}^2)}^2
 \le C\int_{\mathbb{T}^2}\rho v^2\,dx+C\|\nabla v\|_{L^2(\mathbb{T}^2)}^2.
\end{equation}
\end{lemma}

Under the hypotheses of Lemma~\ref{lem:weighted-poincare}, $\eqref{eq:periodic-interpolation}_2$ and \eqref{eq:weighted-poincare-l2} give, for $2\le m<\infty$,
\begin{equation}\label{eq:weighted-poincare}
 \|v\|_{L^m}\le C\sqrt m\bigl(\|\sqrt{\rho}\,v\|_{L^2}+\|\nabla v\|_{L^2}\bigr),
\end{equation}
where $C$ is independent of $m$.

We use the following interpolation inequality on $\mathbb{T}^2$
due to Fan, Jiang, Li, and L\"u~\cite{bib6}.

\begin{lemma}\label{lem:uniform-interpolation}
Let $p,m\in[8,\infty)$ and let $f\in C^\infty(\mathbb{T}^2;\mathbb{R})$
satisfy $\int_{\mathbb{T}^2} f\,dx=0$. Then
\begin{equation}\label{eq:exponent-interpolation}
  \|f\|_{L^\infty}
  \le C\sqrt{p}\,\|f\|_{L^2}^{\eta}
     \|\nabla f\|_{L^2}^{1-\eta-\sigma}
     \|\nabla f\|_{L^m}^{\sigma},
\end{equation}
where $C>0$ is independent of $f$, $p$, and $m$, and
\[
  \eta=\frac{2(m-2)}{p(m-2)+2m},
  \qquad
  \sigma=\frac{2m}{p(m-2)+2m}.
\]
\end{lemma}

We write
\begin{align}
 \mathrm{BMO}(\mathbb{T}^2)
 &=\bigl\{f\in L^1_{\mathrm{loc}}(\mathbb{T}^2)
          :\|f\|_{\mathrm{BMO}}<\infty\bigr\},
 \label{eq:periodic-bmo}\\
 \|f\|_{\mathrm{BMO}}
 &=\sup_{\substack{x\in\mathbb{T}^2\\0<r<d}}
   \frac{1}{|\Omega_r(x)|}
   \int_{\Omega_r(x)}
   \left|f(y)-\frac{1}{|\Omega_r(x)|}
                 \int_{\Omega_r(x)}f(z)\,dz\right|\,dy.
 \label{eq:periodic-bmo-seminorm}
\end{align}
Here $d$ is the diameter of $\mathbb{T}^2$, and $\Omega_r(x)$ is the
periodic ball of center $x$ and radius $r$.
Next, let \(\Delta^{-1}\) denote the Laplacian inverse with zero mean on \(\mathbb T^2\) and 
\begin{align}
    R_i=(-\Delta)^{-1/2}\partial_i\label{eq:riesz-transform}
\end{align}
be the usual Riesz transform on \(\mathbb T^2\), satisfy
\begin{align}
 R_iR_j&=\partial_i\partial_j(-\Delta)^{-1}=R_jR_i=(R_iR_j)^* .
 \label{eq:riesz-symmetry}
\end{align}
For a linear operator $L$, write
\begin{equation}\label{eq:commutator}
 [b,L](f)=bLf-L(bf).
\end{equation}

The following estimate follows by polarization from
\cite{PradoliniSalinas2007}.

\begin{lemma}\label{lem:double-commutator}
Let $b,c\in C^\infty(\mathbb{T}^2;\mathbb{R})$. Then there exists a
constant $C$ depending only on $\mathbb{T}^2$ such that, for every
$f\in L^2(\mathbb{T}^2)$,
\begin{equation}\label{eq:iterated-commutator}
 \|[b,[c,R_iR_j]]f\|_{L^2}
 \le C\|b\|_{\mathrm{BMO}}\|c\|_{\mathrm{BMO}}\|f\|_{L^2}.
\end{equation}
\end{lemma}
\begin{proof}
For sufficiently
small $\varepsilon>0$, the John--Nirenberg inequality gives a uniform bound on the Muckenhoupt characteristic of order two
for the weight $\exp(2\operatorname{Re}(zb+wc))$ whenever
\[
 |z|\|b\|_{\mathrm{BMO}}+|w|\|c\|_{\mathrm{BMO}}\le\varepsilon .
\]
The weighted Calder\'on--Zygmund estimate on $\mathbb{T}^2$
then bounds $e^{zb+wc}R_iR_je^{-zb-wc}$ uniformly on $L^2$ in this range.
Cauchy's integral formula in $z$ and $w$ at $(0,0)$ yields \eqref{eq:iterated-commutator}.
See~\cite{bib10} for the Euclidean argument.
\end{proof}

For the paired estimate, set
\begin{equation}\label{eq:flux-commutator}
 F=\sum_{i,j=1}^2[u_i,R_iR_j](\rho u_j).
\end{equation}
\begin{lemma}\label{lem:paired-commutator}
Let $\rho, \mathbf{u}\in C^\infty(\mathbb{T}^2)$.
Then
\begin{align}
 \int\rho F\,dx
 &=-\frac12\sum_{i,j=1}^2\int\rho[u_i,[u_j,R_iR_j]]\rho\,dx,
 \label{eq:paired-commutator-identity}\\
 \left|\int\rho F\,dx\right|
 &\le C\|\nabla\mathbf{u}\|_{L^2}^2\|\rho\|_{L^2}^2.
 \label{eq:paired-commutator-estimate}
\end{align}
\end{lemma}

\begin{proof}
By \eqref{eq:riesz-symmetry} and \eqref{eq:flux-commutator},
\[
\begin{aligned}
 \sum_{i,j=1}^2\int\rho[u_i,[u_j,R_iR_j]]\rho\,dx
 &=2\sum_{i,j=1}^2\int
   \bigl(\rho R_iR_j(\rho u_i u_j)-\rho u_iR_iR_j(\rho u_j)\bigr)\,dx\\
 &=-2\int\rho F\,dx.
\end{aligned}
\]
By \eqref{eq:paired-commutator-identity}, \eqref{eq:iterated-commutator},
and Poincar\'e's inequality ,
\[
 \left|\int\rho F\,dx\right|
 \le C\|\rho\|_{L^2}^2\sum_{i,j=1}^2
       \|u_i\|_{\mathrm{BMO}}\|u_j\|_{\mathrm{BMO}}
 \le C\|\nabla\mathbf{u}\|_{L^2}^2\|\rho\|_{L^2}^2.
\]
\end{proof}
The following Beale--Kato--Majda type inequality was proved in
\cite{Kato1986} for $\operatorname{div}\mathbf{u}=0$ (see also \cite{BealeKatoMajda1984}).
\begin{lemma}[See \cite{BealeKatoMajda1984,Kato1986,HuangLiXin2011}]\label{lem:log-div-curl}
For $2<q<\infty$, there is a constant $C(q)$ such that the following estimate holds
for all $\nabla \mathbf{u}\in W^{1,q}(\mathbb{T}^2)$,
\begin{equation}\label{eq:log-div-curl}
 \|\nabla \mathbf{u}\|_{L^\infty}
 \le C\bigl(\|\operatorname{div}\mathbf{u}\|_{L^\infty}
           +\|\operatorname{rot}\mathbf{u}\|_{L^\infty}\bigr)
       \log\!\bigl(e+\|\nabla^2\mathbf{u}\|_{L^q}\bigr)
       +C\|\nabla \mathbf{u}\|_{L^2}+C.
\end{equation}
\end{lemma}
\section{A priori estimates (I): an upper bound for the density}\label{sec:3}
In this section and the next, we will always assume that $(\rho_0,m_0)$
satisfies \eqref{eq:local-existence-data} and $(\rho,\mathbf{u})$ is the
strong solution to \eqref{eq:system}--\eqref{eq:initial-data} on
$\mathbb{T}^2\times(0,T_1]$ obtained by Lemma~\ref{lem:local-existence}.
Here $0<T_1\le T$ and $T$ is fixed.

We use the effective viscous flux, vorticity,
and energy quantities of Huang and Li~\cite{bib3}. Under
\eqref{eq:pressure-law-viscosity} we have
\begin{equation}\label{eq:linear-coefficient}
 2\mu+\lambda(\rho)=2\mu+\rho,\qquad \nabla^\perp=(\partial_2,-\partial_1).
\end{equation}
The effective viscous flux and vorticity are
\begin{equation}\label{eq:flux-vorticity}
 G=(2\mu+\lambda(\rho))\operatorname{div}\mathbf{u}-(P-\bar P),\qquad
 \omega=\nabla^\perp\cdot \mathbf{u}=\partial_2 u_1-\partial_1 u_2 .
\end{equation}
We define
\begin{align}
 A_1(t)^2&=\int\left(\omega^2+\frac{G^2}{2\mu+\lambda(\rho)}\right)\,dx,
  \label{eq:flux-norm}\\
 A_2(t)^2&=\int\rho|\dot{\mathbf{u}}|^2\,dx,\label{eq:acceleration-norm}\\
 A_3(t)^2&=\int\big((2\mu+\lambda(\rho))(\operatorname{div}\mathbf{u})^2+\mu\omega^2\big)\,dx.
  \label{eq:velocity-dissipation}
\end{align}
The pressure estimates also use
\begin{equation}\label{eq:unshifted-flux}
 \widetilde G=G-\bar P=(2\mu+\lambda(\rho))\operatorname{div}\mathbf{u}-P .
\end{equation}
Set
\begin{equation}\label{eq:density-maximum}
 R(t)=\max\!\left\{e^2,\ \sup_{0\le s\le t}\|\rho(\cdot,s)\|_{L^\infty}\right\}.
\end{equation}
The viscosity primitive at $\beta=1$ satisfies
\begin{equation}\label{eq:viscosity-primitive}
 \theta(s)=2\mu\log s+s,\qquad s\theta'(s)=2\mu+s,\qquad s>0.
\end{equation}
\eqref{eq:continuity} and \eqref{eq:unshifted-flux} therefore give
\begin{equation}\label{eq:primitive-transport}
 \frac{D}{Dt}\theta(\rho)+P=-\widetilde G .
\end{equation}
Throughout the estimates, $C$ depends on $\mu$, $\gamma$,
$\mathbb{T}^2$, $M$, $\|\rho_0\|_{L^\infty}$, and
$\|\mathbf{u}_0\|_{H^1}$. The constant $C_T$ may also depend on $T$
and is uniform for all $T_1\le T$.

Integration of \eqref{eq:system} over $\mathbb{T}^2$ gives conservation
of mass and momentum. We may assume that
\begin{equation}\label{eq:zero-momentum}
  \int \rho\mathbf{u}\,dx=0.
\end{equation}
\begin{lemma}\label{lem:basic-energy}
There exists a constant $C$ such that
\begin{equation}\label{eq:basic-energy}
 \sup_{0\le t\le T_1}\int\left(\tfrac12\rho|\mathbf{u}|^2+\frac{P}{\gamma-1}\right)\,dx
 +\int_0^{T_1} A_3(t)^2\,dt\le C .
\end{equation}
\end{lemma}
\begin{proof}
The periodic div--curl identity reads
\begin{equation}\label{eq:div-curl}
 \|\nabla \mathbf{u}\|_{L^2}^2=\|\operatorname{div}\mathbf{u}\|_{L^2}^2+\|\omega\|_{L^2}^2 .
\end{equation}
Testing \eqref{eq:momentum} by $\mathbf{u}$ and using \eqref{eq:continuity} and \eqref{eq:div-curl} give
\begin{equation}\label{eq:energy-balance}
 \frac{d}{dt}\int\left(\tfrac12\rho|\mathbf{u}|^2+\frac{P}{\gamma-1}\right)\,dx
 +A_3(t)^2=0 .
\end{equation}
Integration in time yields \eqref{eq:basic-energy}.
\end{proof}
By \eqref{eq:basic-energy}, we have
\begin{equation}\label{eq:mean-flux}
 |\overline{\widetilde G}|=\left|\int(\rho\operatorname{div}\mathbf{u}-P)\right|
 \le M^{1/2}\left(\int\rho(\operatorname{div}\mathbf{u})^2\right)^{1/2}
      +\int P\le C(1+A_3(t)) .
\end{equation}
\eqref{eq:basic-energy} and~\eqref{eq:zero-momentum} yield
\[
  M\overline{\mathbf{u}}
  =-\int\rho(\mathbf{u}-\overline{\mathbf{u}})\,dx.
\]
H\"older's inequality and the Poincar\'e--Sobolev inequality yield
\[
  |\overline{\mathbf{u}}|
  \le \frac{\|\rho\|_{L^\gamma}}{M}
      \|\mathbf{u}-\overline{\mathbf{u}}\|_{L^{\gamma/(\gamma-1)}}
  \le C\|\nabla\mathbf{u}\|_{L^2}.
\]
In particular,
\[
  \|\mathbf{u}\|_{L^2}
  \le \|\mathbf{u}-\overline{\mathbf{u}}\|_{L^2}
      +|\mathbb{T}^2|^{1/2}|\overline{\mathbf{u}}|
  \le C\|\nabla\mathbf{u}\|_{L^2}.
\]
 \eqref{eq:div-curl} and 
\eqref{eq:velocity-dissipation} give
\[
  \mu\|\nabla\mathbf{u}\|_{L^2}^2
  =\mu\bigl(\|\operatorname{div}\mathbf{u}\|_{L^2}^2
            +\|\omega\|_{L^2}^2\bigr)
  \le A_3(t)^2.
\]
Thus
\begin{equation}\label{eq:velocity-h1}
  \|\mathbf{u}\|_{H^1}\le C A_3(t).
\end{equation}

Define the momentum potential by
\begin{equation}\label{eq:momentum-potential}
 \psi=(-\Delta)^{-1}\operatorname{div}(\rho \mathbf{u}).
\end{equation}
By \eqref{eq:riesz-symmetry}, the commutator in \eqref{eq:flux-commutator} satisfies
\begin{equation}\label{eq:potential-commutator}
 F=\mathbf{u}\cdot\nabla\psi-(-\Delta)^{-1}\operatorname{div}\operatorname{div}(\rho \mathbf{u}\otimes \mathbf{u}).
\end{equation}
\begin{proposition}\label{prop:quadratic-density}
For $\gamma>1$, there exists a constant $C_T>0$ such that
\begin{equation}\label{eq:quadratic-density}
 \sup_{0\le t\le T_1}\|\rho(t)\|_{L^2}^2
 +\int_0^{T_1}\!\int \rho^{\gamma+1}\,dx\,dt
 \le C_T,
\end{equation}
where $C_T$ is independent of $R(T_1)$ and of the positive lower bound for $\rho$.
\end{proposition}

\begin{proof}
Apply $(-\Delta)^{-1}\operatorname{div}$ to \eqref{eq:momentum} and use
\eqref{eq:potential-commutator} to obtain
\begin{equation}\label{eq:potential-transport}
  \frac{D}{Dt}\psi=-\widetilde G+\overline{\widetilde G}+F.
\end{equation}
Subtract \eqref{eq:potential-transport} from \eqref{eq:primitive-transport}
and multiply the resulting equation by $\rho$. Integrating  and using
\eqref{eq:continuity}, we have
\begin{equation}\label{eq:density-moment}
  \frac{d}{dt}\int \rho\bigl(\theta(\rho)-\psi\bigr)\,dx
  +\int \rho^{\gamma+1}\,dx
  =-M\overline{\widetilde G}-\int \rho F\,dx
  =:J_1+J_2.
\end{equation}
\eqref{eq:mean-flux} gives
\[
  |J_1|\le C\bigl(1+A_3(t)\bigr).
\]
For $J_2$, the paired commutator bound \eqref{eq:paired-commutator-estimate} and \eqref{eq:velocity-h1} give
\[
  |J_2|
  \le C\|\nabla\mathbf{u}\|_{L^2}^2\|\rho\|_{L^2}^2
  \le C A_3(t)^2\|\rho\|_{L^2}^2.
\]
The elliptic estimate for \eqref{eq:momentum-potential}, together with
\eqref{eq:basic-energy}, gives
\begin{equation}\label{eq:momentum-potential-l4}
  \|\psi\|_{L^4}
  \le C\|\rho\mathbf{u}\|_{L^{4/3}}
  \le C\|\sqrt{\rho}\,\mathbf{u}\|_{L^2}\|\rho\|_{L^2}^{1/2}
  \le C\|\rho\|_{L^2}^{1/2}.
\end{equation}
From \eqref{eq:momentum-potential-l4}, we obtain
\[
  \left|\int \rho\psi\,dx\right|
  \le \|\rho\|_{L^{4/3}}\|\psi\|_{L^4}\le M^{1/2}\|\rho\|_{L^2}^{1/2}\|\psi\|_{L^4}
  \le C\|\rho\|_{L^2}.
\]
By \eqref{eq:viscosity-primitive}, \eqref{eq:basic-energy}, and
\eqref{eq:momentum-potential-l4},
\[
 \left|\int \rho\bigl(\theta(\rho)-\psi\bigr)\,dx
       -\|\rho\|_{L^2}^2\right|
 \le C\bigl(1+\|\rho\|_{L^2}\bigr).
\]
Young's inequality allows us to choose a fixed $c_\rho>0$ such that
\begin{equation}\label{eq:moment-comparison}
  c\bigl(1+\|\rho\|_{L^2}^2\bigr)
  \le I_\rho(t):=c_\rho+\int \rho\bigl(\theta(\rho)-\psi\bigr)\,dx
  \le C\bigl(1+\|\rho\|_{L^2}^2\bigr).
\end{equation}
Insert the bounds for $J_1$ and $J_2$ into \eqref{eq:density-moment}.
Together with \eqref{eq:moment-comparison}, this gives
\[
  I_\rho'(t)+\int \rho^{\gamma+1}\,dx
  \le C\bigl(1+A_3(t)^2\bigr)I_\rho(t).
\]
Gronwall's inequality and \eqref{eq:basic-energy}
then give \eqref{eq:quadratic-density}.
\end{proof}
Define
\begin{equation}\label{eq:pressure-energy}
 J(t):=1+\int\frac{P^2}{2\mu+\lambda(\rho)}\,dx .
\end{equation}
Its integrand, as a function of the density, is
\begin{equation}\label{eq:pressure-integrand}
 j(s)=\frac{s^{2\gamma}}{2\mu+s}.
\end{equation}
For $\gamma>2$, we use a correction to the pressure integral. In the two ranges of $\gamma$ the energy is
\begin{equation}\label{eq:corrected-energy}
\begin{gathered}
 \mathcal{V}=e+\int\left(\mu\omega^2+\frac{\widetilde G^2}{2\mu+\lambda(\rho)}\right)\,dx
 +\begin{cases}
 J-1,&1<\gamma\le2,\\[2pt]
 \displaystyle\int j_*(\rho)\,dx,&\gamma>2,
 \end{cases}\\
 j_*(r)=r\int_0^r
 \left(\gamma-1-\frac{s}{2\mu+s}\right)\frac{2s^{2\gamma-2}}{2\mu+s}\,ds .
\end{gathered}
\end{equation}
For $\gamma>2$ and $0<s\le r$, we have
\[
\begin{gathered}
 2(\gamma-2)\le 2\left(\gamma-1-\frac{s}{2\mu+s}\right)\le 2(\gamma-1),\\
 \frac1{2\mu+r}\le \frac1{2\mu+s}\le\frac{r}{s(2\mu+r)} .
\end{gathered}
\]
Integrating against $s^{2\gamma-2}$ gives $$cj(r)\le j_*(r)\le Cj(r).$$ Hence \eqref{eq:unshifted-flux} gives, 
\begin{equation}\label{eq:energy-comparison}
 c(J+A_3(t)^2)\le \mathcal{V}\le C(J+A_3(t)^2),\qquad \mathcal{V}\ge e.
\end{equation}
By \eqref{eq:unshifted-flux}, \eqref{eq:basic-energy}, and the
definitions \eqref{eq:flux-norm} and \eqref{eq:velocity-dissipation}, there are constants $c,C>0$ such that
\begin{equation}\label{eq:energy-flux-comparison}
 c(e+A_1(t)^2+A_3(t)^2)\le \mathcal{V}\le C(e+A_1(t)^2+A_3(t)^2).
\end{equation}

\begin{lemma}\label{lem:elliptic-flux}
For $G$ and $\omega$ in \eqref{eq:flux-vorticity}, there exist constants
$C_p$ and $C$ such that
\begin{equation}\label{eq:elliptic-flux}
\begin{aligned}
 \|\nabla G\|_{L^p}+\|\nabla\omega\|_{L^p}
 &\le C_p\|\rho\dot{\mathbf{u}}\|_{L^p},\qquad 1<p<\infty,\\
 \|\nabla G\|_{L^2}+\|\nabla\omega\|_{L^2}
 &\le C\sqrt R A_2(t).
\end{aligned}
\end{equation}
\end{lemma}
\begin{proof}
The momentum equation takes the form
\begin{equation}\label{eq:momentum-flux}
 \rho\dot{\mathbf{u}}=\nabla G+\mu\nabla^\perp\omega .
\end{equation}
Taking its divergence and curl gives
\begin{align}
 \Delta G&=\operatorname{div}(\rho\dot{\mathbf{u}}),\label{eq:flux-poisson}\\
 \mu\Delta\omega&=\nabla^\perp\cdot(\rho\dot{\mathbf{u}}).\label{eq:vorticity-poisson}
\end{align}
Periodic elliptic estimates give the first inequality in \eqref{eq:elliptic-flux}.
 For $p=2$, use
$\|\rho\dot{\mathbf{u}}\|_{L^2}\le\sqrt{R}\,A_2(t)$ to obtain the second inequality.
\end{proof}

Define the coefficient
\begin{equation}\label{eq:flux-energy-coefficient}
 c(s)=1+\frac{s}{2\mu+s},\qquad s\ge0 .
\end{equation}
\begin{lemma}\label{lem:flux-identity}
With $c$ as in \eqref{eq:flux-energy-coefficient}, we have
\begin{equation}\label{eq:flux-energy-identity}
\begin{aligned}
 &\frac{d}{dt}\int\left(\mu\omega^2+\frac{\widetilde G^2}{2\mu+\lambda(\rho)}\right)\,dx+2A_2(t)^2\\
 &=\int c(\rho)\frac{\widetilde G^2}{2\mu+\lambda(\rho)}
      \operatorname{div}\mathbf{u}\,dx
   -2\int c(\rho)\widetilde G(\operatorname{div}\mathbf{u})^2\,dx\\
 &\quad+2\gamma\int\frac{P\widetilde G}{2\mu+\lambda(\rho)}
      \operatorname{div}\mathbf{u}\,dx
   +4\int\widetilde G\det\nabla\mathbf{u}\,dx
   -\mu\int\omega^2\operatorname{div}\mathbf{u}\,dx .
\end{aligned}
\end{equation}
\end{lemma}
\begin{proof}
As in \cite{bib3},
differentiation and \eqref{eq:unshifted-flux} give
\begin{equation}\label{eq:material-div-curl}
\begin{aligned}
 \nabla^\perp\cdot\dot{\mathbf{u}}
 &=\frac{D\omega}{Dt}+\omega\operatorname{div}\mathbf{u},\\
 \operatorname{div}\dot{\mathbf{u}}
 &=\frac{D}{Dt}\left(\frac{\widetilde G}{2\mu+\lambda(\rho)}\right)
   +\frac{D}{Dt}\left(\frac{P}{2\mu+\lambda(\rho)}\right)\\
 &\quad-2\det\nabla\mathbf{u}+(\operatorname{div}\mathbf{u})^2 .
\end{aligned}
\end{equation}
Since $\nabla G=\nabla\widetilde G$, testing \eqref{eq:momentum-flux} against
$2\dot{\mathbf{u}}$ and integrating over $\mathbb{T}^2$ gives
\begin{equation}\label{eq:flux-energy-test}
 2A_2(t)^2
 =-2\int\widetilde G\operatorname{div}\dot{\mathbf{u}}\,dx
  -2\mu\int\omega\nabla^\perp\cdot\dot{\mathbf{u}}\,dx .
\end{equation}
\eqref{eq:continuity} and~\eqref{eq:pressure-law} give
\begin{equation}\label{eq:pressure-quotient-transport}
 \frac{D}{Dt}\left(\frac{P}{2\mu+\lambda(\rho)}\right)
 =\left(\frac{\rho}{2\mu+\lambda(\rho)}-\gamma\right)
  \frac{P}{2\mu+\lambda(\rho)}
  \operatorname{div}\mathbf{u}.
\end{equation}
Substituting \eqref{eq:material-div-curl} and \eqref{eq:pressure-quotient-transport}
into \eqref{eq:flux-energy-test} and integrating by parts,
we obtain
\begin{equation}\label{eq:flux-energy-expanded}
\begin{aligned}
 &\frac{d}{dt}\int\left(\mu\omega^2+
   \frac{\widetilde G^2}{2\mu+\lambda(\rho)}\right)\,dx+2A_2(t)^2\\
 &=-\mu\int\omega^2\operatorname{div}\mathbf{u}\,dx
   +4\int\widetilde G\det\nabla\mathbf{u}\,dx\\
 &\quad-2\int\widetilde G(\operatorname{div}\mathbf{u})^2\,dx
   +2\mu\int\frac{\widetilde G^2}{(2\mu+\lambda(\rho))^2}
     \operatorname{div}\mathbf{u}\,dx\\
 &\quad-2\int\frac{\rho P\widetilde G}{(2\mu+\lambda(\rho))^2}
     \operatorname{div}\mathbf{u}\,dx
   +2\gamma\int\frac{P\widetilde G}{2\mu+\lambda(\rho)}
     \operatorname{div}\mathbf{u}\,dx .
\end{aligned}
\end{equation}
\eqref{eq:unshifted-flux} and \eqref{eq:linear-coefficient} give
\begin{equation}\label{eq:flux-energy-regrouping}
\begin{aligned}
 &\frac{2\mu\widetilde G^2}{(2\mu+\lambda(\rho))^2}\operatorname{div}\mathbf{u}
 -2\widetilde G(\operatorname{div}\mathbf{u})^2
 -\frac{2\rho P\widetilde G}{(2\mu+\lambda(\rho))^2}\operatorname{div}\mathbf{u}\\
 &=c(\rho)\frac{\widetilde G^2}{2\mu+\lambda(\rho)}\operatorname{div}\mathbf{u}
   -2c(\rho)\widetilde G(\operatorname{div}\mathbf{u})^2.
\end{aligned}
\end{equation}
Applying \eqref{eq:flux-energy-regrouping} to \eqref{eq:flux-energy-expanded}
gives \eqref{eq:flux-energy-identity}.
\end{proof}

\begin{lemma}\label{lem:nonlinear-flux}
For every $\varepsilon>0$, there exists a constant $C_\varepsilon$
such that
\begin{equation}\label{eq:flux-div-bound}
 \|\operatorname{div}\mathbf{u}\|_{L^2}\left\|\frac{\widetilde G^2}{2\mu+\lambda(\rho)}\right\|_{L^2}
 \le\varepsilon A_2(t)^2+C_\varepsilon (R\log R) A_3(t)^2 \mathcal{V}+C\mathcal{V}.
\end{equation}
\begin{equation}\label{eq:jacobian-bound}
 \left|\int\widetilde G\det\nabla \mathbf{u}\right|
 \le C\|\nabla G\|_{L^2}\|\nabla \mathbf{u}\|_{L^2}^{2}
 \le\varepsilon A_2(t)^2+C_\varepsilon R A_3(t)^4
 \le\varepsilon A_2(t)^2+C_\varepsilon (R\log R) A_3(t)^2 \mathcal{V}.
\end{equation}
\begin{equation}\label{eq:vorticity-div-bound}
 \left|\int\omega^2\operatorname{div}\mathbf{u}\right|
 \le CA_3(t)\sqrt{\mathcal{V}}\sqrt R A_2(t)\le\varepsilon A_2(t)^2+C_\varepsilon (R\log R) A_3(t)^2 \mathcal{V} .
\end{equation}
\end{lemma}
\begin{proof}
By \eqref{eq:mean-flux} and \eqref{eq:energy-comparison}, we can bound $\overline{\widetilde G}$.
Then \eqref{eq:corrected-energy} gives
\begin{equation}\label{eq:flux-l2}
 |\overline{\widetilde G}|\le C\sqrt{\mathcal{V}},\qquad
 \|\widetilde G\|_{L^2}\le C\sqrt{R\mathcal{V}},\qquad
 \|\omega\|_{L^2}\le C\sqrt{\mathcal{V}}.
\end{equation}
For $0<\delta\le1/4$, H\"older's inequality, \eqref{eq:periodic-interpolation}
and \eqref{eq:flux-l2} give
\begin{equation}\label{eq:flux-product}
\begin{aligned}
 \left\|\frac{\widetilde G^2}{2\mu+\lambda(\rho)}\right\|_{L^2}
 &\le C\mathcal{V}^{(1-\delta)/2}\|\widetilde G\|_{L^{2(1+\delta)/\delta}}^{1+\delta}\\
 &\le C\delta^{-1/2}\mathcal{V}^{(1-\delta)/2}
 \|G-\bar G\|_{L^2}^{\delta}\|\nabla G\|_{L^2}
 +C\mathcal{V}^{(1-\delta)/2}|\overline{\widetilde G}|^{1+\delta}\\
 &\le C\delta^{-1/2}R^{\delta/2}\sqrt{\mathcal{V}}\,\|\nabla G\|_{L^2}+C\mathcal{V}.
\end{aligned}
\end{equation}
The second inequality uses the uniform bound
\[
 \left(\frac{2(1+\delta)}{\delta}\right)^{(1+\delta)/2}\le C\delta^{-1/2}.
\]
Choose $\delta=(4\log R)^{-1}$, so $R^{\delta/2}$ is bounded.
Applying Young's inequality to \eqref{eq:elliptic-flux} gives \eqref{eq:flux-div-bound}.
We have
\[
 CA_3(t)\mathcal{V}\le C(R\log R)A_3(t)^2\mathcal{V}+C\mathcal{V}.
\]
As in \cite{bib3}, we use
\begin{gather*}
 \det\nabla\mathbf{u}
 =\nabla u_1\cdot\nabla^\perp u_2,
 \qquad \int\det\nabla\mathbf{u}\,dx=0,\\
 \|\det\nabla\mathbf{u}\|_{\mathcal H^1}
 \le C\|\nabla\mathbf{u}\|_{L^2}^2.
\end{gather*}
Hardy--BMO duality, \eqref{eq:elliptic-flux}, \eqref{eq:velocity-h1},
\eqref{eq:energy-comparison}, and Young's inequality imply
\begin{align*}
 \left|\int\widetilde G\det\nabla\mathbf{u}\,dx\right|
 &=\left|\int G\det\nabla\mathbf{u}\,dx\right|\\
 &\le C\|G\|_{\mathrm{BMO}}\|\det\nabla\mathbf{u}\|_{\mathcal H^1}
 \le C\|\nabla G\|_{L^2}\|\nabla\mathbf{u}\|_{L^2}^2\\
 &\le C\sqrt R\,A_2(t)A_3(t)^2
 \le \varepsilon A_2(t)^2+C_\varepsilon R A_3(t)^4\\
 &\le \varepsilon A_2(t)^2+C_\varepsilon(R\log R)A_3(t)^2\mathcal V,
\end{align*}
which proves \eqref{eq:jacobian-bound}.

Applying Young's inequality to the bounds from
\eqref{eq:periodic-interpolation}, \eqref{eq:velocity-dissipation},
\eqref{eq:flux-l2}, and \eqref{eq:elliptic-flux}, we obtain
\begin{align*}
 \left|\int\omega^2\operatorname{div}\mathbf{u}\,dx\right|
 &\le \|\omega\|_{L^4}^2\|\operatorname{div}\mathbf{u}\|_{L^2}\\
 &\le C\|\omega\|_{L^2}\|\nabla\omega\|_{L^2}
        \|\operatorname{div}\mathbf{u}\|_{L^2}\\
 &\le C A_3(t)\sqrt{\mathcal V}\sqrt R\,A_2(t)
 \le \varepsilon A_2(t)^2+C_\varepsilon(R\log R)A_3(t)^2\mathcal V .
\end{align*}
\end{proof}
\begin{lemma}\label{lem:pressure-identity}
For $j$ in \eqref{eq:pressure-integrand}, define
\begin{equation}\label{eq:pressure-renormalization}
 a(\rho)=\rho j'(\rho)-j(\rho)=\left(2\gamma-1-\frac{\rho}{2\mu+\lambda(\rho)}\right)j(\rho).
\end{equation}
In particular,
\[
 0\le a\le Cj.
\]
For $J$ in
\eqref{eq:pressure-energy}, we have the identities
\begin{equation}\label{eq:pressure-identities}
\begin{aligned}
 J'&=-\int a\operatorname{div}\mathbf{u},\\
 J'+\int\frac{aP}{2\mu+\lambda(\rho)}&=-\int\frac{a\widetilde G}{2\mu+\lambda(\rho)}.
\end{aligned}
\end{equation}
\end{lemma}
\begin{proof}
The bound on $a$ follows from \eqref{eq:linear-coefficient} and $\gamma>1$.
The continuity equation \eqref{eq:continuity} and periodic integration by
parts give the first identity. The second follows from \eqref{eq:unshifted-flux}.
\end{proof}

\begin{proposition}\label{prop:flux-energy}
For $1<\gamma\le2$, there exist constants $c>0$ and $C_T>0$ such that
\begin{align}
 \int_0^{T_1}\mathcal{V}\,dt &\le C_T,\label{eq:subquadratic-energy-integral}\\
 \mathcal{V}'+cA_2(t)^2 &\le C_T(R\log R)\mathcal{V}^2,
 \label{eq:subquadratic-energy-differential}\\
 \mathcal{V}^2 &\le C_T R^2\bigl(1+A_2(t)^2\bigr).
 \label{eq:subquadratic-energy-pointwise}
\end{align}
\end{proposition}
\begin{proof}
It follows from \eqref{eq:pressure-integrand} and \eqref{eq:quadratic-density} that
\begin{align}
 \int_0^{T_1}(J-1)\,dt
 &\le C\int_0^{T_1}\!\int(1+\rho^{\gamma+1})\,dx\,dt
 \le C_T,\label{eq:pressure-integrability}\\
 J(t)-1
 &\le C\bigl(1+R^{(2\gamma-3)_+}\|\rho(t)\|_{L^2}^2\bigr)
 \le C_T R^{(2\gamma-3)_+}.\label{eq:pressure-pointwise}
\end{align}
The estimate \eqref{eq:subquadratic-energy-integral} follows from
\eqref{eq:pressure-integrability}, \eqref{eq:basic-energy}, and
\eqref{eq:energy-comparison}. Testing \eqref{eq:momentum-flux} against
$\mathbf{u}$ and using \eqref{eq:basic-energy}, we have
\begin{equation}\label{eq:dissipation-test}
\begin{aligned}
 A_3(t)^2
 &=\int P\operatorname{div}\mathbf{u}\,dx
   -\int \rho\dot{\mathbf{u}}\cdot\mathbf{u}\,dx\\
 &\le A_3(t)\sqrt{J-1}+C A_2(t)
 \le \tfrac12 A_3(t)^2+C\bigl(J-1+A_2(t)\bigr).
\end{aligned}
\end{equation}
For \eqref{eq:subquadratic-energy-pointwise}, use \eqref{eq:energy-comparison},
\eqref{eq:dissipation-test}, and \eqref{eq:pressure-pointwise} to obtain
\begin{equation}\label{eq:subquadratic-energy-from-pressure}
\begin{aligned}
 \mathcal{V}^2
 &\le C\bigl(1+(J-1)^2+A_2(t)^2\bigr)\\
 &\le C_T R^{2(2\gamma-3)_+}\bigl(1+A_2(t)^2\bigr)
 \le C_T R^2\bigl(1+A_2(t)^2\bigr).
\end{aligned}
\end{equation}
We derive \eqref{eq:subquadratic-energy-differential} from
\begin{equation}\label{eq:pressure-growth-weight}
 \frac{j(\rho)^2}{2\mu+\lambda(\rho)}
 =j(\rho)\left(\frac{P}{2\mu+\lambda(\rho)}\right)^2
 \le R^{2\gamma-2}j(\rho).
\end{equation}
From \eqref{eq:pressure-identities}, Lemma~\ref{lem:pressure-identity},
\eqref{eq:pressure-growth-weight}, and \eqref{eq:corrected-energy}, we obtain
\begin{equation}\label{eq:pressure-growth}
\begin{aligned}
 J'
 &\le C\left\|\frac{j(\rho)}{\sqrt{2\mu+\lambda(\rho)}}\right\|_{L^2}
        \left\|\frac{\widetilde G}{\sqrt{2\mu+\lambda(\rho)}}\right\|_{L^2}\\
 &\le C R^{\gamma-1}\sqrt{J-1}
        \left\|\frac{\widetilde G}{\sqrt{2\mu+\lambda(\rho)}}\right\|_{L^2}
 \le C R\mathcal{V}.
\end{aligned}
\end{equation}
Using \eqref{eq:unshifted-flux}, we arrive at
\begin{equation}\label{eq:subquadratic-divergence-bound}
\begin{aligned}
 &C\int\frac{\widetilde G^2+P|\widetilde G|}{2\mu+\lambda(\rho)}
            |\operatorname{div}\mathbf{u}|\,dx\le C\|\operatorname{div}\mathbf{u}\|_{L^2}
 \left(
  \left\|\frac{\widetilde G^2}{2\mu+\lambda(\rho)}\right\|_{L^2}
  +\left\|\frac{\widetilde G P}{2\mu+\lambda(\rho)}\right\|_{L^2}
 \right).
\end{aligned}
\end{equation}
For $1<\gamma<2$, take
\begin{equation}\label{eq:subcritical-pressure-exponents}
 r=\frac{2}{\gamma-1},\qquad q=\frac{2}{2-\gamma},\qquad
 \frac1r+\frac1q=\frac12.
\end{equation}
From \eqref{eq:pressure-law} and \eqref{eq:quadratic-density},
\begin{equation}\label{eq:subcritical-pressure-factor}
 \left\|\frac{P}{2\mu+\lambda(\rho)}\right\|_{L^r}
 \le \|\rho\|_{L^2}^{\gamma-1}
 \le C_T.
\end{equation}
Using \eqref{eq:periodic-interpolation}, \eqref{eq:elliptic-flux}, and \eqref{eq:flux-l2},
we have
\begin{equation}\label{eq:subcritical-flux-factor}
\begin{aligned}
 \|\widetilde G\|_{L^q}
 &\le C_T\|G-\bar G\|_{L^2}^{2-\gamma}\|\nabla G\|_{L^2}^{\gamma-1}
      +|\overline{\widetilde G}|\\
 &\le C_T R^{1/2}\mathcal{V}^{(2-\gamma)/2}A_2(t)^{\gamma-1}
      +C_T\sqrt{\mathcal{V}}.
\end{aligned}
\end{equation}
H\"older's and Young's inequalities, \eqref{eq:energy-comparison},
and \eqref{eq:subcritical-pressure-exponents}--\eqref{eq:subcritical-flux-factor}
imply
\begin{equation}\label{eq:subcritical-mixed-pressure}
\begin{aligned}
 \|\operatorname{div}\mathbf{u}\|_{L^2}
 \left\|\frac{\widetilde G P}{2\mu+\lambda(\rho)}\right\|_{L^2}
 &\le C_T R^{1/2}\mathcal{V}^{(3-\gamma)/2}A_2(t)^{\gamma-1}
      +C_T\mathcal{V}\\
 &\le \varepsilon A_2(t)^2
      +C_{T,\varepsilon}R^{1/(3-\gamma)}\mathcal{V}+C_T\mathcal{V}\\
 &\le \varepsilon A_2(t)^2+C_{T,\varepsilon}(R\log R)\mathcal{V}^2.
\end{aligned}
\end{equation}
For $\gamma=2$, take
\begin{equation}\label{eq:endpoint-pressure-exponents}
 \delta=(4\log R)^{-1},\qquad
 r=2+\delta,\qquad q=\frac{2(2+\delta)}{\delta},\qquad
 \frac1r+\frac1q=\frac12.
\end{equation}
From \eqref{eq:pressure-law}, \eqref{eq:density-maximum}, and \eqref{eq:quadratic-density},
\begin{equation}\label{eq:endpoint-pressure-factor}
 \left\|\frac{P}{2\mu+\lambda(\rho)}\right\|_{L^r}
 \le\|\rho\|_{L^r}
 \le R^{\delta/(2+\delta)}\|\rho\|_{L^2}^{2/(2+\delta)}
 \le C_T R^{\delta/(2+\delta)}.
\end{equation}
Applying \eqref{eq:periodic-interpolation},
\eqref{eq:elliptic-flux}, and \eqref{eq:flux-l2} to obtain
\begin{equation}\label{eq:endpoint-flux-factor}
\begin{aligned}
 \|\widetilde G\|_{L^q}
 &\le C\sqrt q\,\|G-\bar G\|_{L^2}^{2/q}
                \|\nabla G\|_{L^2}^{1-2/q}+|\overline{\widetilde G}|\\
 &\le C_T\delta^{-1/2}R^{1/2}\mathcal{V}^{\delta/[2(2+\delta)]}
           A_2(t)^{2/(2+\delta)}+C_T\sqrt{\mathcal{V}}.
\end{aligned}
\end{equation}
H\"older's and Young's inequalities with
\eqref{eq:endpoint-pressure-exponents}--\eqref{eq:endpoint-flux-factor}
and \eqref{eq:energy-comparison} imply
\begin{equation}\label{eq:endpoint-mixed-pressure}
\begin{aligned}
 \|\operatorname{div}\mathbf{u}\|_{L^2}
 \left\|\frac{\widetilde G P}{2\mu+\lambda(\rho)}\right\|_{L^2}
 &\le \varepsilon A_2(t)^2+C_{T,\varepsilon}\delta^{-1}
       R^{1+\delta/[2(1+\delta)]}\mathcal{V}
       +C_T R^{\delta/(2+\delta)}\mathcal{V}\\
 &\le \varepsilon A_2(t)^2+C_{T,\varepsilon}(R\log R)\mathcal{V}^2.
\end{aligned}
\end{equation}
In \eqref{eq:flux-energy-identity}, use \eqref{eq:subquadratic-divergence-bound},
\eqref{eq:flux-div-bound}--\eqref{eq:vorticity-div-bound}, and
\eqref{eq:subcritical-mixed-pressure} or \eqref{eq:endpoint-mixed-pressure}
as appropriate. Adding \eqref{eq:pressure-growth} and taking $\varepsilon$
sufficiently small proves \eqref{eq:subquadratic-energy-differential}.
\end{proof}

\begin{proposition}\label{prop:corrected-energy}
For $\gamma>2$, there exist constants $c>0$ and $C>0$ such that
\begin{align}
 \mathcal{V}'+cA_2(t)^2+c\int\frac{P^3}{(2\mu+\lambda(\rho))^2}\,dx
 &\le C(R\log R)A_3(t)^2\mathcal{V}+C\mathcal{V},
 \label{eq:corrected-dissipation}\\
 (\log\mathcal{V})'+c\frac{A_2(t)^2}{\mathcal{V}}+c\mathcal{V}
 &\le C(R\log R)A_3(t)^2+C,
 \label{eq:logarithmic-energy}\\
 \mathcal{V}^2 &\le CR^{2\gamma-2}\bigl(1+A_2(t)^2\bigr).
 \label{eq:superquadratic-pointwise}
\end{align}
\end{proposition}
\begin{proof}
From \eqref{eq:unshifted-flux},
\begin{equation}\label{eq:pressure-expansion}
\begin{aligned}
 &c(\rho)\frac{\widetilde G^2}{2\mu+\lambda(\rho)}
       \operatorname{div}\mathbf{u}
 -2c(\rho)\widetilde G(\operatorname{div}\mathbf{u})^2
 +2\gamma\frac{P\widetilde G}{2\mu+\lambda(\rho)}
       \operatorname{div}\mathbf{u}\\
 &\quad=-c(\rho)\frac{\widetilde G^3}{(2\mu+\lambda(\rho))^2}
 +(2\gamma-3c(\rho))\frac{\widetilde G^2P}{(2\mu+\lambda(\rho))^2}
 +2(\gamma-c(\rho))\frac{\widetilde GP^2}{(2\mu+\lambda(\rho))^2}.
\end{aligned}
\end{equation}
Differentiating the definition of $j_*$ in \eqref{eq:corrected-energy},
we have
\begin{equation}\label{eq:pressure-correction-identity}
 \rho j_*'(\rho)-j_*(\rho)
 =2\left(\gamma-1-\frac{\rho}{2\mu+\rho}\right)
      \frac{\rho^{2\gamma}}{2\mu+\rho}
 =2(\gamma-c(\rho))\frac{P^2}{2\mu+\lambda(\rho)}.
\end{equation}
From \eqref{eq:continuity} and \eqref{eq:pressure-correction-identity},
\begin{equation}\label{eq:pressure-correction-transport}
\begin{aligned}
 \frac{d}{dt}\int j_*(\rho)\,dx
 &=-\int\bigl(\rho j_*'(\rho)-j_*(\rho)\bigr)
              \operatorname{div}\mathbf{u}\,dx\\
 &=-2\int(\gamma-c(\rho))\frac{P^2}{2\mu+\lambda(\rho)}
              \operatorname{div}\mathbf{u}\,dx.
\end{aligned}
\end{equation}
Substitute \eqref{eq:unshifted-flux} into
\eqref{eq:pressure-correction-transport} to obtain
\begin{equation}\label{eq:correction-cancellation}
\begin{aligned}
 \frac{d}{dt}\int j_*(\rho)\,dx
 +2\int(\gamma-c(\rho))\frac{P^3}{(2\mu+\lambda(\rho))^2}\,dx
 &=-2\int(\gamma-c(\rho))\frac{\widetilde GP^2}{(2\mu+\lambda(\rho))^2}\,dx .
\end{aligned}
\end{equation}
Adding \eqref{eq:correction-cancellation} to \eqref{eq:flux-energy-identity}
and using \eqref{eq:pressure-expansion}, we obtain
\begin{equation}\label{eq:corrected-flux-identity}
\begin{aligned}
 &\mathcal{V}'+2A_2(t)^2
   +2\int(\gamma-c(\rho))\frac{P^3}{(2\mu+\lambda(\rho))^2}\,dx\\
 &\quad=-\int c(\rho)\frac{\widetilde G^3}{(2\mu+\lambda(\rho))^2}\,dx
   +\int(2\gamma-3c(\rho))\frac{\widetilde G^2P}{(2\mu+\lambda(\rho))^2}\,dx\\
 &\qquad+4\int\widetilde G\det\nabla\mathbf{u}\,dx
   -\mu\int\omega^2\operatorname{div}\mathbf{u}\,dx .
\end{aligned}
\end{equation}
The coefficients satisfy
\begin{equation}\label{eq:pressure-dissipation-coefficient}
 1\le c(\rho)\le 2,\qquad
 2(\gamma-c(\rho))\ge 2(\gamma-2)>0.
\end{equation}
H\"older's inequality and
\eqref{eq:basic-energy}, \eqref{eq:corrected-energy}, and
\eqref{eq:elliptic-flux} imply
\begin{align}
 \|\widetilde G\|_{L^{2\gamma/(\gamma+1)}}
 &\le
 \left\|\frac{\widetilde G}{\sqrt{2\mu+\lambda(\rho)}}\right\|_{L^2}
 \left\|\sqrt{2\mu+\lambda(\rho)}\right\|_{L^{2\gamma}}
 \le C\sqrt{\mathcal{V}},
 \label{eq:flux-pstar}\\
 \|\nabla G\|_{L^{2\gamma/(\gamma+1)}}
 &\le C\|\rho\dot{\mathbf{u}}\|_{L^{2\gamma/(\gamma+1)}}
 \le C\|\sqrt{\rho}\|_{L^{2\gamma}}A_2(t)
 \le CA_2(t).
 \label{eq:flux-gradient-pstar}
\end{align}
The Sobolev embedding and
\eqref{eq:flux-pstar}--\eqref{eq:flux-gradient-pstar} yield
\begin{equation}\label{eq:flux-l2gamma}
 \|\widetilde G\|_{L^{2\gamma}}
 \le C\bigl(\|\widetilde G\|_{L^{2\gamma/(\gamma+1)}}
           +\|\nabla G\|_{L^{2\gamma/(\gamma+1)}}\bigr)
 \le C\bigl(\sqrt{\mathcal{V}}+A_2(t)\bigr).
\end{equation}
H\"older's inequality, Young's inequality, with
\eqref{eq:pressure-law}, \eqref{eq:basic-energy}, \eqref{eq:flux-pstar}, and
\eqref{eq:flux-l2gamma}, imply
\begin{equation}\label{eq:quadratic-flux-pressure}
\begin{aligned}
 \int\frac{\widetilde G^2P}{(2\mu+\lambda(\rho))^2}\,dx
 &\le \int \widetilde G^2\rho^{\gamma-2}\,dx
 \le \|\widetilde G\|_{L^\gamma}^2\|\rho\|_{L^\gamma}^{\gamma-2}\\
 &\le C\|\widetilde G\|_{L^{2\gamma/(\gamma+1)}}^{2/\gamma}
        \|\widetilde G\|_{L^{2\gamma}}^{2(\gamma-1)/\gamma}\\
 &\le C\mathcal{V}^{1/\gamma}A_2(t)^{2(\gamma-1)/\gamma}+C\mathcal{V}\\
 &\le \varepsilon A_2(t)^2+C_\varepsilon\mathcal{V}.
\end{aligned}
\end{equation}
From \eqref{eq:unshifted-flux},
\begin{equation}\label{eq:cubic-flux-decomposition}
 \frac{\widetilde G^3}{(2\mu+\lambda(\rho))^2}
 =\frac{\widetilde G^2}{2\mu+\lambda(\rho)}\operatorname{div}\mathbf{u}
 -\frac{\widetilde G^2P}{(2\mu+\lambda(\rho))^2}.
\end{equation}
It follows from \eqref{eq:cubic-flux-decomposition},
\eqref{eq:pressure-dissipation-coefficient}, \eqref{eq:flux-div-bound}, and
\eqref{eq:quadratic-flux-pressure} that
\begin{equation}\label{eq:cubic-flux-bound}
\begin{aligned}
 \left|\int c(\rho)\frac{\widetilde G^3}{(2\mu+\lambda(\rho))^2}\,dx\right|
 &\le C\|\operatorname{div}\mathbf{u}\|_{L^2}
       \left\|\frac{\widetilde G^2}{2\mu+\lambda(\rho)}\right\|_{L^2}
      +C\int\frac{\widetilde G^2P}{(2\mu+\lambda(\rho))^2}\,dx\\
 &\le \varepsilon A_2(t)^2
      +C_\varepsilon(R\log R)A_3(t)^2\mathcal{V}
      +C_\varepsilon\mathcal{V}.
\end{aligned}
\end{equation}
Use \eqref{eq:quadratic-flux-pressure}, \eqref{eq:cubic-flux-bound}, and
\eqref{eq:jacobian-bound}--\eqref{eq:vorticity-div-bound} in
\eqref{eq:corrected-flux-identity}. Taking $\varepsilon$ sufficiently small
and applying \eqref{eq:pressure-dissipation-coefficient} proves
\eqref{eq:corrected-dissipation}.

Using \eqref{eq:pressure-energy} and \eqref{eq:basic-energy}, we obtain
\begin{equation}\label{eq:pressure-moment-dissipation}
 (J-1)^2\le\left(\int P\,dx\right)
            \int\frac{P^3}{(2\mu+\lambda(\rho))^2}\,dx
 \le C\int\frac{P^3}{(2\mu+\lambda(\rho))^2}\,dx.
\end{equation}
From \eqref{eq:energy-comparison} and \eqref{eq:pressure-moment-dissipation},
\begin{equation}\label{eq:energy-from-pressure-dissipation}
\begin{aligned}
 \mathcal{V}^2
 &\le C\bigl((J-1)^2+(1+A_3(t)^2)^2\bigr)\\
 &\le C\int\frac{P^3}{(2\mu+\lambda(\rho))^2}\,dx
       +C\mathcal{V}\bigl(1+A_3(t)^2\bigr).
\end{aligned}
\end{equation}
Divide \eqref{eq:corrected-dissipation} by $\mathcal{V}$ and use
\eqref{eq:energy-from-pressure-dissipation} to estimate
\begin{equation}\label{eq:logarithmic-energy-calculation}
\begin{aligned}
 (\log\mathcal{V})'+c\frac{A_2(t)^2}{\mathcal{V}}+c\mathcal{V}
 &\le C(R\log R)A_3(t)^2+C\bigl(1+A_3(t)^2\bigr)\\
 &\le C(R\log R)A_3(t)^2+C.
\end{aligned}
\end{equation}
This is \eqref{eq:logarithmic-energy}.

From \eqref{eq:pressure-energy}, \eqref{eq:density-maximum}, and
\eqref{eq:basic-energy}, we have
\begin{equation}\label{eq:superquadratic-pressure-pointwise}
 J-1
 =\int\frac{P^2}{2\mu+\lambda(\rho)}\,dx
 \le R^{\gamma-1}\int P\,dx
 \le CR^{\gamma-1}.
\end{equation}
With \eqref{eq:superquadratic-pressure-pointwise} and \eqref{eq:energy-comparison},
it implies
\begin{equation}\label{eq:superquadratic-pointwise-calculation}
 \mathcal{V}^2
 \le C\bigl(1+(J-1)^2+A_2(t)^2\bigr)
 \le CR^{2\gamma-2}\bigl(1+A_2(t)^2\bigr),
\end{equation}
as claimed in \eqref{eq:superquadratic-pointwise}.
\end{proof}

\begin{proposition}\label{prop:material-derivative}
For $\gamma>1$. There exist constants $c>0$, $C_T>0$, and
$\alpha_\gamma\ge1$, with $\alpha_\gamma$ depending only on $\gamma$,
such that
\begin{equation}\label{eq:acceleration-inequality}
 (A_2(t)^2)'+c\|\nabla\dot{\mathbf{u}}\|_{L^2}^2
 \le C_T R^{\alpha_\gamma}\mathcal{V}\bigl(1+A_2(t)^2\bigr).
\end{equation}
For $\gamma>2$, $C_T$ depends only on the fixed parameters and
the bound in \eqref{eq:basic-energy}.
\end{proposition}
\begin{proof}
Following \cite{bib3}, apply
$\partial_t+\operatorname{div}(\mathbf{u}\,\cdot)$ to the material form
of \eqref{eq:momentum} and test the $j$th component with $\dot u_j$.
From \eqref{eq:continuity},
\begin{equation}\label{eq:acceleration-test}
\begin{aligned}
 \tfrac12(A_2(t)^2)'
 &=-\int\dot u_j
 \bigl[\partial_j P_t+\operatorname{div}(\mathbf{u}\,\partial_j P)\bigr]\,dx\\
 &\quad+\mu\int\dot u_j
 \bigl[\partial_t\Delta u_j+\operatorname{div}(\mathbf{u}\,\Delta u_j)\bigr]\,dx\\
 &\quad+\int\dot u_j\Bigl[
 \partial_j\partial_t\bigl((\mu+\lambda(\rho))\operatorname{div}\mathbf{u}\bigr)\\
 &\qquad
 +\operatorname{div}\bigl(\mathbf{u}\,\partial_j
                ((\mu+\lambda(\rho))\operatorname{div}\mathbf{u})\bigr)
 \Bigr]\,dx\\
 &=:N_1+N_2+N_3.
\end{aligned}
\end{equation}
It follows from \eqref{eq:continuity} and \eqref{eq:pressure-law} that
\begin{equation}\label{eq:acceleration-pressure-transport}
 P_t+\operatorname{div}(P\mathbf{u})
 =-(\gamma-1)P\operatorname{div}\mathbf{u}.
\end{equation}
For $N_1$ in \eqref{eq:acceleration-test}, use
\eqref{eq:acceleration-pressure-transport} and integrate by parts to obtain
\begin{equation}\label{eq:acceleration-pressure-bound}
\begin{aligned}
 N_1
 &=-\int P\Bigl((\gamma-1)\operatorname{div}\mathbf{u}\,
             \operatorname{div}\dot{\mathbf{u}}
              +\partial_j u_k\,\partial_k\dot u_j\Bigr)\,dx\\
 &\le C R^\gamma\|\nabla\mathbf{u}\|_{L^2}
                     \|\nabla\dot{\mathbf{u}}\|_{L^2}\\
 &\le \tfrac{\mu}{8}\|\nabla\dot{\mathbf{u}}\|_{L^2}^2
        +CR^{2\gamma}\|\nabla\mathbf{u}\|_{L^2}^2 .
\end{aligned}
\end{equation}
Integration by parts in $N_2$ gives
\begin{equation}\label{eq:acceleration-shear-bound}
\begin{aligned}
 N_2
 &=-\mu\int|\nabla\dot{\mathbf{u}}|^2\,dx
   -\mu\int\partial_i\dot u_j
       \bigl(\partial_i u_j\,\operatorname{div}\mathbf{u}
                  -\partial_i u_k\,\partial_k u_j\bigr)\,dx\\
 &\quad+\mu\int\partial_k\dot u_j\,\partial_i u_k\,\partial_i u_j\,dx\\
 &\le -\tfrac{3\mu}{4}\|\nabla\dot{\mathbf{u}}\|_{L^2}^2
         +C\|\nabla\mathbf{u}\|_{L^4}^4 .
\end{aligned}
\end{equation}
From \eqref{eq:continuity}, \eqref{eq:linear-coefficient}, and
\eqref{eq:material-div-curl},
\begin{equation}\label{eq:acceleration-bulk-transport}
\begin{aligned}
 &\partial_t\bigl((\mu+\lambda(\rho))\operatorname{div}\mathbf{u}\bigr)
 +\operatorname{div}\bigl(\mathbf{u}(\mu+\lambda(\rho))
                                 \operatorname{div}\mathbf{u}\bigr)\\
 &\quad=(\mu+\lambda(\rho))
       \bigl(\operatorname{div}\dot{\mathbf{u}}-\partial_j u_k\,\partial_k u_j\bigr)
      +\mu(\operatorname{div}\mathbf{u})^2 .
\end{aligned}
\end{equation}
Integrating $N_3$ by parts and inserting
\eqref{eq:acceleration-bulk-transport}, we find
\begin{equation}\label{eq:acceleration-bulk-bound}
\begin{aligned}
 N_3
 &=-\int(\mu+\lambda(\rho))(\operatorname{div}\dot{\mathbf{u}})^2\,dx\\
 &\quad+\int\operatorname{div}\dot{\mathbf{u}}\Bigl[
       (\mu+\lambda(\rho))\partial_j u_k\,\partial_k u_j
       -\mu(\operatorname{div}\mathbf{u})^2\Bigr]\,dx\\
 &\quad+\int(\mu+\lambda(\rho))\operatorname{div}\mathbf{u}\,
                   \partial_j u_k\,\partial_k\dot u_j\,dx\\
 &\le -\int(\mu+\lambda(\rho))(\operatorname{div}\dot{\mathbf{u}})^2\,dx
       +\tfrac{\mu}{8}\|\nabla\dot{\mathbf{u}}\|_{L^2}^2
       +CR^2\|\nabla\mathbf{u}\|_{L^4}^4 .
\end{aligned}
\end{equation}
Insert \eqref{eq:acceleration-pressure-bound}, \eqref{eq:acceleration-shear-bound},
and \eqref{eq:acceleration-bulk-bound} into \eqref{eq:acceleration-test} to obtain
\begin{equation}\label{eq:acceleration-dissipation}
\begin{aligned}
 &(A_2(t)^2)' +\mu\|\nabla\dot{\mathbf{u}}\|_{L^2}^2
 +2\int(\mu+\lambda(\rho))(\operatorname{div}\dot{\mathbf{u}})^2\,dx\\
 &\qquad\le CR^2\|\nabla\mathbf{u}\|_{L^4}^4
        +CR^{2\gamma}\|\nabla\mathbf{u}\|_{L^2}^2.
\end{aligned}
\end{equation}
The div--curl estimate and \eqref{eq:unshifted-flux} imply
\begin{equation}\label{eq:acceleration-velocity-l4}
\begin{aligned}
 \|\nabla\mathbf{u}\|_{L^4}^4
 &\le C\bigl(\|\operatorname{div}\mathbf{u}\|_{L^4}^4+\|\omega\|_{L^4}^4\bigr)\\
 &\le C\bigl(\|\widetilde G\|_{L^4}^4+\|\omega\|_{L^4}^4+R^{4\gamma}\bigr).
\end{aligned}
\end{equation}
Using \eqref{eq:periodic-interpolation}, \eqref{eq:elliptic-flux},
and \eqref{eq:flux-l2}, we have
\begin{equation}\label{eq:acceleration-flux-l4}
\begin{aligned}
 \|\widetilde G\|_{L^4}^4+\|\omega\|_{L^4}^4
 &\le C\Bigl(\|G-\bar G\|_{L^2}^2\|\nabla G\|_{L^2}^2
       +\|\omega\|_{L^2}^2\|\nabla\omega\|_{L^2}^2
       +|\overline{\widetilde G}|^4\Bigr)\\
 &\le CR^2\mathcal{V}A_2(t)^2+C\mathcal{V}^2 .
\end{aligned}
\end{equation}
Use \eqref{eq:acceleration-velocity-l4}, \eqref{eq:acceleration-flux-l4},
\eqref{eq:div-curl}, and \eqref{eq:energy-comparison} in
\eqref{eq:acceleration-dissipation}. Then \eqref{eq:subquadratic-energy-pointwise}
and \eqref{eq:superquadratic-pointwise} imply
\begin{equation}\label{eq:acceleration-energy-bound}
\begin{aligned}
 (A_2(t)^2)' +\mu\|\nabla\dot{\mathbf{u}}\|_{L^2}^2
 &\le C\Bigl(R^4\mathcal{V}A_2(t)^2+R^2\mathcal{V}^2
                 +R^{2+4\gamma}+R^{2\gamma}\mathcal{V}\Bigr)\\
 &\le C_T R^{2+4\gamma}\mathcal{V}\bigl(1+A_2(t)^2\bigr).
\end{aligned}
\end{equation}
Take $\alpha_\gamma=2+4\gamma$ in \eqref{eq:acceleration-inequality}.
\end{proof}
\begin{lemma}\label{lem:weighted-energy}
Let $R(T_1)\le\rho_\sharp$ with $\rho_\sharp\ge e^2$. Suppose the
following inequalities hold on $[0,T_1]$,
\begin{align}
 \mathcal{V}'+cA_2(s)^2
 &\le C_{\rho_\sharp}\mathcal{V}^2,\label{eq:bounded-density-ode}\\
 (A_2(s)^2)'+c\|\nabla\dot{\mathbf{u}}\|_{L^2}^2
 &\le C_{\rho_\sharp}\mathcal{V}\bigl(1+A_2(s)^2\bigr).
 \label{eq:bounded-density-acceleration-ode}
\end{align}
Then there exists a constant $C_{\rho_\sharp}$, which may depend on $T$
as well, such that
\begin{align}
 \sup_{0\le s\le T_1}
 \bigl(\mathcal{V}(s)+\|\mathbf{u}(s)\|_{H^1}^2\bigr)
 &\le C_{\rho_\sharp},\label{eq:bounded-density-estimates}\\
 \int_0^{T_1}A_2(s)^2\,ds
 &\le C_{\rho_\sharp},\label{eq:bounded-density-acceleration-integral}\\
 \sup_{0<s\le T_1}s\bigl(1+A_2(s)^2\bigr)
 &\le C_{\rho_\sharp},\label{eq:bounded-density-weighted-acceleration}\\
 \int_0^{T_1}s\|\nabla\dot{\mathbf{u}}\|_{L^2}^2\,ds
 &\le C_{\rho_\sharp}.\label{eq:bounded-density-weighted-gradient}
\end{align}
\end{lemma}

\begin{proof}
The initial data and \eqref{eq:energy-comparison} imply
\begin{equation}\label{eq:bounded-density-initial-energy}
 \mathcal{V}(0)\le C\bigl(J(0)+A_3(0)^2\bigr)\le C.
\end{equation}
From \eqref{eq:pressure-energy}, \eqref{eq:energy-comparison}, and
\eqref{eq:basic-energy},
\begin{equation}\label{eq:bounded-density-energy-integral}
 \int_0^{T_1}\mathcal{V}\,ds
 \le C\int_0^{T_1}\bigl(J+A_3(s)^2\bigr)\,ds
 \le C_{\rho_\sharp}.
\end{equation}
Integrating \eqref{eq:bounded-density-ode} after division by $\mathcal{V}$,
\begin{equation}\label{eq:bounded-density-log-energy}
 \log\mathcal{V}(s)
 \le \log\mathcal{V}(0)+C_{\rho_\sharp}\int_0^s\mathcal{V}(\xi)\,d\xi
 \le C_{\rho_\sharp}.
\end{equation}
\eqref{eq:velocity-h1} and \eqref{eq:energy-comparison} imply
\eqref{eq:bounded-density-estimates}. Integrating \eqref{eq:bounded-density-ode},
\begin{equation}\label{eq:bounded-density-acceleration-energy}
 c\int_0^{T_1}A_2(s)^2\,ds
 \le \mathcal{V}(0)
   +C_{\rho_\sharp}\sup_{0\le s\le T_1}\mathcal{V}(s)
                  \int_0^{T_1}\mathcal{V}\,ds
 \le C_{\rho_\sharp}.
\end{equation}
Multiplying \eqref{eq:bounded-density-acceleration-ode} by $s$, we obtain
\begin{equation}\label{eq:weighted-acceleration}
 \frac{d}{ds}\bigl[s(1+A_2(s)^2)\bigr]+cs\|\nabla\dot{\mathbf{u}}\|_{L^2}^2
 \le 1+A_2(s)^2+C_{\rho_\sharp}\mathcal{V}\,s\bigl(1+A_2(s)^2\bigr).
\end{equation}
The weighted quantity vanishes at $s=0$. Gronwall's inequality,
\eqref{eq:bounded-density-energy-integral}, and
\eqref{eq:bounded-density-acceleration-integral} imply
\eqref{eq:bounded-density-weighted-acceleration}. Integration of
\eqref{eq:weighted-acceleration} proves \eqref{eq:bounded-density-weighted-gradient}.
\end{proof}

If $R(T_1)\le\rho_\sharp$, the hypotheses of Lemma~\ref{lem:weighted-energy}
follow from \eqref{eq:subquadratic-energy-differential} or
\eqref{eq:corrected-dissipation}, and \eqref{eq:acceleration-inequality}.

\begin{proposition}\label{prop:initial-interval}
There exist $\tau_0\in(0,1]$ and $C_0>0$, depending only on $\mu$,
$\gamma$, $\mathbb{T}^2$, and the initial data, such that
\begin{equation}\label{eq:initial-interval}
 \sup_{0\le s\le\min\{T_1,\tau_0\}}
 \bigl[R(s)+\mathcal{V}(s)+s(1+A_2(s)^2)\bigr]\le C_0.
\end{equation}
If $T_1\ge\tau_0/2$, then
\begin{equation}\label{eq:initial-interval-start}
 R(\tau_0/2)+\mathcal{V}(\tau_0/2)+A_2(\tau_0/2)\le C_0.
\end{equation}
\end{proposition}

\begin{proof}
 Set
\begin{equation}\label{eq:initial-density-threshold}
 \rho_\sharp=2\max\bigl\{e^2,\,1+\|\rho_0\|_{L^\infty}\bigr\}.
\end{equation}
Use $T=1$ in Lemma~\ref{lem:weighted-energy}. Let $t\le\min\{T_1,1\}$ and suppose $R(t)\le\rho_\sharp$. On $[0,t]$,
\eqref{eq:bounded-density-estimates}, \eqref{eq:flux-l2},
\eqref{eq:elliptic-flux}, and \eqref{eq:weighted-poincare} imply
\begin{align}
 \|\widetilde G\|_{L^2}&\le C_{\rho_\sharp},
 \label{eq:initial-flux-l2}\\
 \|\nabla G\|_{L^4}
 &\le C_{\rho_\sharp}\bigl(A_2(s)+\|\nabla\dot{\mathbf{u}}\|_{L^2}\bigr).
 \label{eq:initial-flux-gradient}
\end{align}
The Gagliardo--Nirenberg inequality and \eqref{eq:initial-flux-l2}--\eqref{eq:initial-flux-gradient} yield
\begin{equation}\label{eq:initial-flux-pointwise}
\begin{aligned}
 \|\widetilde G\|_{L^\infty}
 &\le C\|\widetilde G\|_{L^2}^{1/3}\|\nabla G\|_{L^4}^{2/3}
       +C\|\widetilde G\|_{L^2}\\
 &\le C_{\rho_\sharp}\bigl(1+A_2(s)^{2/3}
                     +\|\nabla\dot{\mathbf{u}}\|_{L^2}^{2/3}\bigr).
\end{aligned}
\end{equation}
From \eqref{eq:initial-flux-pointwise},
\eqref{eq:bounded-density-acceleration-integral}, and
\eqref{eq:bounded-density-weighted-gradient},
\begin{equation}\label{eq:initial-flux-integral}
\begin{aligned}
 \int_0^t\|\widetilde G\|_{L^\infty}\,ds
 &\le C_{\rho_\sharp}\Bigg[
 t+t^{2/3}\left(\int_0^t A_2(s)^2\,ds\right)^{1/3}\\
 &\qquad\qquad
 +\left(\int_0^t s\|\nabla\dot{\mathbf{u}}\|_{L^2}^2\,ds\right)^{1/3}
  \left(\int_0^t s^{-1/2}\,ds\right)^{2/3}\Bigg]\\
 &\le C_{\rho_\sharp}t^{1/3}.
\end{aligned}
\end{equation}
At a first contact with $\rho_\sharp$, integrating
\eqref{eq:primitive-transport} along a trajectory gives
\begin{equation}\label{eq:initial-density-gap}
 \theta(\rho_\sharp)-\theta\bigl(1+\|\rho_0\|_{L^\infty}\bigr)
 \le \int_0^t\|\widetilde G\|_{L^\infty}\,ds.
\end{equation}
Choose $\tau_0\in(0,1]$ such that
\begin{equation}\label{eq:initial-time-choice}
 C_{\rho_\sharp}\tau_0^{1/3}
 <\theta(\rho_\sharp)-\theta\bigl(1+\|\rho_0\|_{L^\infty}\bigr).
\end{equation}
\eqref{eq:initial-flux-integral}--\eqref{eq:initial-time-choice} exclude contact
on $[0,\min\{T_1,\tau_0\}]$. The bounds \eqref{eq:bounded-density-estimates}
and \eqref{eq:bounded-density-weighted-acceleration} prove
\eqref{eq:initial-interval}. Enlarge $C_0$ for \eqref{eq:initial-interval-start}.
\end{proof}

For $T_1>\tau_0$, use \eqref{eq:initial-interval-start} at
\begin{equation}\label{eq:density-start-time}
 \tau=\frac{\tau_0}{2}.
\end{equation}

\begin{proposition}\label{prop:weighted-flux}
Suppose the bounds in \eqref{eq:initial-interval-start} hold at $\tau>0$.
With $\alpha_\gamma$ as in \eqref{eq:acceleration-inequality}, there exists
$C_{\tau,T}>0$, independent of $R(T_1)$ and the positive lower bound for
$\rho$, such that for $\tau\le T_1\le T$,
\begin{align}
 \int_\tau^{T_1}\frac{\mathcal{V}+e+A_2(t)^2}{\mathcal{V} R\log R}\,dt
 &\le C_{\tau,T},\label{eq:weighted-energies}\\
 \int_\tau^{T_1}\frac{1}{R^{\alpha_\gamma+2}}
 \left(1+\frac{\|\nabla\dot{\mathbf{u}}\|_{L^2}^2}{e+A_2(t)^2}\right)\,dt
 &\le C_{\tau,T},\label{eq:weighted-gradient-integral}\\
 \log\!\left(1+\frac{e+A_2(t)^2}{\mathcal{V}(t)}\right)
 &\le C_{\tau,T}R(t)^{\alpha_\gamma+2}
 \qquad(\tau\le t\le T_1),\label{eq:weighted-acceleration-log}\\
 \int_\tau^{T_1}\frac{\|\widetilde G\|_{L^\infty}^2}
                       {\mathcal{V}(R\log R)^2}\,dt
 &\le C_{\tau,T}.\label{eq:weighted-flux}
\end{align}
\end{proposition}

\begin{proof}
The running maximum $R$ is locally absolutely continuous and nondecreasing.
For $1<\gamma\le2$, divide \eqref{eq:subquadratic-energy-differential}
by $\mathcal{V}R\log R$ to obtain
\begin{equation}\label{eq:weighted-log-subquadratic}
\begin{aligned}
 \left(\frac{\log\mathcal{V}}{R\log R}\right)'
 +\frac{cA_2(t)^2}{\mathcal{V}R\log R}
 &\le C_T\mathcal{V}
     -\frac{\log\mathcal{V}}{(R\log R)^2}(R\log R)'\\
 &\le C_T\mathcal{V}.
\end{aligned}
\end{equation}
For $\gamma>2$, division of \eqref{eq:logarithmic-energy} by $R\log R$ yields
\begin{equation}\label{eq:weighted-log-superquadratic}
 \left(\frac{\log\mathcal{V}}{R\log R}\right)'
 +c\left(\frac{A_2(t)^2}{\mathcal{V}R\log R}
         +\frac{\mathcal{V}}{R\log R}\right)
 \le CA_3(t)^2+\frac{C}{R\log R}.
\end{equation}
Integrating \eqref{eq:weighted-log-subquadratic} and
\eqref{eq:weighted-log-superquadratic} with
\eqref{eq:subquadratic-energy-integral} and \eqref{eq:basic-energy},
respectively, yields
\begin{equation}\label{eq:weighted-log-integral}
 \int_\tau^{T_1}\left(
    \frac{A_2(t)^2}{\mathcal{V}R\log R}+\frac{\mathcal{V}}{R\log R}
 \right)\,dt\le C_{\tau,T}.
\end{equation}
\eqref{eq:energy-comparison} and \eqref{eq:weighted-log-integral} imply
\eqref{eq:weighted-energies}.

For \eqref{eq:weighted-gradient-integral}, divide
\eqref{eq:acceleration-inequality} by $(e+A_2(t)^2)R^{\alpha_\gamma+2}$ to obtain
\begin{equation}\label{eq:weighted-gradient-differential}
\begin{aligned}
 \left(\frac{\log(e+A_2(t)^2)}{R^{\alpha_\gamma+2}}\right)'
 +\frac{c\|\nabla\dot{\mathbf{u}}\|_{L^2}^2}
        {(e+A_2(t)^2)R^{\alpha_\gamma+2}}
 &\le C_T\frac{\mathcal{V}}{R^2}
 \le C_T\frac{\mathcal{V}}{R\log R}.
\end{aligned}
\end{equation}
Integrate \eqref{eq:weighted-gradient-differential} and use
\eqref{eq:weighted-log-integral}. For \eqref{eq:weighted-acceleration-log},
divide \eqref{eq:acceleration-inequality} by $e+A_2(t)^2$ and integrate to obtain
\begin{equation}\label{eq:acceleration-log-integrated}
\begin{aligned}
 \log(e+A_2(t)^2)
 &\le \log(e+A_2(\tau)^2)
     +C_T R(t)^{\alpha_\gamma+1}\log R(t)
                 \int_\tau^t\frac{\mathcal{V}}{R\log R}\,ds\\
 &\le C_{\tau,T}R(t)^{\alpha_\gamma+2}.
\end{aligned}
\end{equation}
\eqref{eq:weighted-acceleration-log} follows from
\eqref{eq:acceleration-log-integrated} and \eqref{eq:energy-comparison}.

For \eqref{eq:weighted-flux}, write
\begin{equation}\label{eq:weighted-ratios}
 X=1+\frac{e+A_2(t)^2}{\mathcal{V}},\qquad
 Y=1+\frac{\|\nabla\dot{\mathbf{u}}\|_{L^2}^2}{e+A_2(t)^2}.
\end{equation}
For $m\ge8$, \eqref{eq:elliptic-flux} and \eqref{eq:weighted-poincare}
give
\begin{equation}\label{eq:weighted-flux-gradient}
 \|\nabla G\|_{L^m}
 \le CmR\|\dot{\mathbf{u}}\|_{L^m}
 \le Cm^{3/2}R\bigl(A_2(t)+\|\nabla\dot{\mathbf{u}}\|_{L^2}\bigr).
\end{equation}
Lemma~\ref{lem:uniform-interpolation}, \eqref{eq:flux-l2},
\eqref{eq:elliptic-flux}, and \eqref{eq:weighted-flux-gradient} imply
\begin{equation}\label{eq:flux-interpolation}
\begin{aligned}
 \frac{\|G-\bar G\|_{L^\infty}^2}{\mathcal{V}}
 &\le Cp\frac{(R\mathcal{V})^\eta
       [R(e+A_2(t)^2)]^{1-\eta-\sigma}
       [m^3 R^2(e+A_2(t)^2)Y]^\sigma}{\mathcal{V}}\\
 &\le Cp R^{1+\sigma}m^{3\sigma}X^{1-\eta}Y^\sigma .
\end{aligned}
\end{equation}
Choose
\begin{equation}\label{eq:weighted-interpolation-exponents}
 p=8+(\alpha_\gamma+7)\log R,\qquad m=e^p .
\end{equation}
\begin{align}
 &(\alpha_\gamma+7)\log R\le p\le(\alpha_\gamma+11)\log R,
 \label{eq:weighted-p-range}\\
 &R^{(\alpha_\gamma+1)\sigma}+m^{3\sigma}\le C,
 \label{eq:weighted-exponent-bounds}\\
 &X^{\sigma-\eta}
 =\exp\!\left(\frac{4\log X}{p(m-2)+2m}\right)
 \le\exp\!\left(\frac{C_{\tau,T}R^{\alpha_\gamma+2}}{pm}\right)
 \le C_{\tau,T}.
 \label{eq:weighted-exponent-difference}
\end{align}
Using \eqref{eq:flux-l2}, \eqref{eq:flux-interpolation}, and
\eqref{eq:weighted-p-range}--\eqref{eq:weighted-exponent-difference},
\begin{equation}\label{eq:weighted-flux-pointwise}
\begin{aligned}
 \frac{\|\widetilde G\|_{L^\infty}^2}{\mathcal{V}(R\log R)^2}
 &\le C_{\tau,T}\frac{X^{1-\sigma}Y^\sigma}{R\log R}+C\\
 &=C_{\tau,T}R^{(\alpha_\gamma+1)\sigma}(\log R)^{-\sigma}
     \left(\frac{X}{R\log R}\right)^{1-\sigma}
     \left(\frac{Y}{R^{\alpha_\gamma+2}}\right)^\sigma+C\\
 &\le C_{\tau,T}\left(\frac{X}{R\log R}
                 +\frac{Y}{R^{\alpha_\gamma+2}}+1\right).
\end{aligned}
\end{equation}
Integrating \eqref{eq:weighted-flux-pointwise} and using
\eqref{eq:weighted-energies} and \eqref{eq:weighted-gradient-integral}
proves \eqref{eq:weighted-flux}.
\end{proof}
\begin{lemma}\label{lem:record-maxima}
Let $(\rho,\mathbf{u})$ be a smooth positive solution on $[0,T_1]$, and
let $R$ be defined by \eqref{eq:density-maximum}. Then
\begin{equation}\label{eq:record-maximum}
 \theta'(R)R'+R^\gamma\le \|\widetilde G\|_{L^\infty}
 \qquad\text{a.e. on }\{R'>0\}.
\end{equation}
Consequently, for almost every $t\in[0,T_1]$,
\begin{align}
 R'&\le \|\widetilde G\|_{L^\infty},
 \label{eq:maximum-derivatives}\\
 4R'R^\gamma&\le \|\widetilde G\|_{L^\infty}^2.
 \label{eq:maximum-derivative-square}
\end{align}
\end{lemma}

\begin{proof}
Let $\mathcal{X}$ be the flow of $\mathbf{u}$ and set
\begin{equation}\label{eq:density-maximum-along-flow}
 Z(t,a)=\theta\bigl(\rho(\mathcal{X}(t,a),t)\bigr),\qquad
 Z_*(t)=\max_{a\in\mathbb{T}^2}Z(t,a)
       =\theta\!\left(\max_{x\in\mathbb{T}^2}\rho(x,t)\right).
\end{equation}
The smoothness of $\rho$ and \eqref{eq:density-maximum} imply that
$Z_*$ and $R$ are Lipschitz on $[0,T_1]$. In particular, $R$ is
absolutely continuous.

If $Z_*(t)<\theta(R(t))$, then $R$ is constant immediately to the
right of $t$. Thus, at a differentiability time with $R'(t)>0$,
a maximizer $a_t$ of $Z(t,\cdot)$ satisfies
\begin{equation}\label{eq:density-maximum-envelope}
 Z(s,a_t)\le Z_*(s)\le\theta(R(s)),\qquad
 Z(t,a_t)=Z_*(t)=\theta(R(t)).
\end{equation}
The difference $\theta(R(s))-Z(s,a_t)$ has a minimum at $s=t$.
Using \eqref{eq:primitive-transport}, we obtain
\begin{equation}\label{eq:density-maximum-contact}
\begin{aligned}
 \theta'(R(t))R'(t)+R(t)^\gamma
 &=\partial_t Z(t,a_t)+R(t)^\gamma\\
 &=-\widetilde G(\mathcal{X}(t,a_t),t)
 \le \|\widetilde G(t)\|_{L^\infty}.
\end{aligned}
\end{equation}
On $\{R'>0\}$, \eqref{eq:density-maximum-contact} and
\eqref{eq:viscosity-primitive} yield
\begin{align}
 R'&\le R'+R^\gamma
 \le\theta'(R)R'+R^\gamma\le\|\widetilde G\|_{L^\infty},
 \label{eq:density-maximum-linear-comparison}\\
 4R'R^\gamma&\le(R'+R^\gamma)^2\le\|\widetilde G\|_{L^\infty}^2.
 \label{eq:density-maximum-square-comparison}
\end{align}
The estimates \eqref{eq:maximum-derivatives} and
\eqref{eq:maximum-derivative-square} also hold on $\{R'=0\}$.
\end{proof}
\subsection{The density bound for \texorpdfstring{$1<\gamma\le2$}{1 < gamma <= 2}}
\label{subsec:subquadratic-density}

\begin{proposition}\label{prop:subquadratic-density}
Suppose $1<\gamma\le2$ and the initial bounds of Proposition~\ref{prop:weighted-flux}
hold at $\tau$. There exists a constant $C_{\tau,T}>0$ such that
\begin{equation}\label{eq:subquadratic-maximum}
 R(T_1)\le C_{\tau,T}\qquad(\tau\le T_1\le T).
\end{equation}
\end{proposition}
\begin{proof}
From \eqref{eq:maximum-derivatives}, \eqref{eq:weighted-flux}, and
\eqref{eq:subquadratic-energy-integral},
\begin{equation}\label{eq:subquadratic-logarithmic-maximum}
\begin{aligned}
 \log\log R(T_1)-\log\log R(\tau)
 &=\int_\tau^{T_1}\frac{R'}{R\log R}\,dt
 \le\int_\tau^{T_1}\frac{\|\widetilde G\|_{L^\infty}}{R\log R}\,dt\\
 &\le\left(\int_\tau^{T_1}\frac{\|\widetilde G\|_{L^\infty}^2}
                         {\mathcal{V}(R\log R)^2}\,dt\right)^{1/2}
      \left(\int_\tau^{T_1}\mathcal{V}\,dt\right)^{1/2}
 \le C_{\tau,T}.
\end{aligned}
\end{equation}
The bound for $R(\tau)$ and \eqref{eq:subquadratic-logarithmic-maximum}
imply \eqref{eq:subquadratic-maximum}.
\end{proof}

\subsection{The density bound for \texorpdfstring{$\gamma>2$}{gamma > 2}}
\label{subsec:superquadratic-density}

\begin{lemma}\label{lem:pressure-variation}
For $\gamma>2$, take
\begin{equation}\label{eq:pressure-variation-exponents}
 k=\gamma-1,\qquad r=\frac{1}{2k}\in(0,\tfrac12).
\end{equation}
For $J$ in \eqref{eq:pressure-energy}, there exists a constant $C$ such that
\begin{align}
 |J'|&\le C R^k A_3(t)\sqrt J,\label{eq:pressure-variation}\\
 (J')_+&\le C A_3(t)\bigl(1+A_2(t)+A_3(t)\bigr)J^{1/2-r}.
 \label{eq:pressure-positive-variation}
\end{align}
\end{lemma}
\begin{proof}
From \eqref{eq:pressure-identities}, \eqref{eq:pressure-renormalization}, and
\eqref{eq:pressure-energy}, we have
\begin{equation}\label{eq:pressure-total-variation-bound}
\begin{aligned}
 |J'|
 &\le A_3(t)\left(\int\frac{a^2}{2\mu+\lambda(\rho)}\,dx\right)^{1/2}\\
 &\le C A_3(t)\left(\int j(\rho)
          \left(\frac{P}{2\mu+\lambda(\rho)}\right)^2dx\right)^{1/2}
 \le C R^k A_3(t)\sqrt J.
\end{aligned}
\end{equation}
For \eqref{eq:pressure-positive-variation}, use the compression region
\begin{equation}\label{eq:compression-region}
 \Omega_-(t)=\{x\in\mathbb{T}^2:\operatorname{div}\mathbf{u}(x,t)<0\}.
\end{equation}
From \eqref{eq:compression-region}, \eqref{eq:unshifted-flux}, and
\eqref{eq:pressure-renormalization},
\begin{equation}\label{eq:compression-integrand}
 \rho^\gamma\le|\widetilde G|
 \quad\text{on }\Omega_-(t),\qquad
 \frac{a^2}{2\mu+\lambda(\rho)}
 \le C\frac{\rho^{4\gamma}}{(2\mu+\lambda(\rho))^3}
 \le C\rho^{4\gamma-3}.
\end{equation}
\eqref{eq:pressure-energy} controls
\begin{equation}\label{eq:compression-density-moment}
 \int\rho^{2\gamma-1}\,dx
 \le C\left(1+\int j(\rho)\,dx\right)=CJ.
\end{equation}
H\"older's inequality, \eqref{eq:compression-integrand}, and
\eqref{eq:compression-density-moment} imply
\begin{equation}\label{eq:compression-weight-bound}
\begin{aligned}
 \int_{\Omega_-}\frac{a^2}{2\mu+\lambda(\rho)}\,dx
 &\le C\int_{\Omega_-}|\widetilde G|^{q_*}
             \rho^{4\gamma-3-\gamma q_*}\,dx\\
 &\le C\|\widetilde G\|_{L^{2\gamma}}^{q_*}
       \left(\int\rho^{2\gamma-1}\,dx\right)^{1-q_*/(2\gamma)}\\
 &\le C\|\widetilde G\|_{L^{2\gamma}}^{q_*}J^{1-q_*/(2\gamma)},
\end{aligned}
\end{equation}
where
\begin{equation}\label{eq:compression-exponent}
 q_*=\frac{4\gamma k}{2\gamma k+1}<2.
\end{equation}
The Poincar\'e--Sobolev inequality with
\eqref{eq:flux-gradient-pstar} and \eqref{eq:mean-flux} implies
\begin{equation}\label{eq:compression-flux-l2gamma}
\begin{aligned}
 \|\widetilde G\|_{L^{2\gamma}}
 &\le\|G-\bar G\|_{L^{2\gamma}}+|\overline{\widetilde G}|\\
 &\le C\|\nabla G\|_{L^{2\gamma/(\gamma+1)}}+|\overline{\widetilde G}|
 \le C\bigl(1+A_2(t)+A_3(t)\bigr).
\end{aligned}
\end{equation}
Applying \eqref{eq:pressure-identities},
\eqref{eq:compression-weight-bound}, and \eqref{eq:compression-flux-l2gamma},
we obtain
\begin{equation}\label{eq:compression-bound}
\begin{aligned}
 (J')_+
 &\le\int_{\Omega_-}a|\operatorname{div}\mathbf{u}|\,dx
 \le A_3(t)\left(\int_{\Omega_-}\frac{a^2}{2\mu+\lambda(\rho)}\,dx\right)^{1/2}\\
 &\le C A_3(t)\bigl(1+A_2(t)+A_3(t)\bigr)^{q_*/2}J^{1/2-q_*/(4\gamma)}.
\end{aligned}
\end{equation}
For the damping term in \eqref{eq:pressure-identities}, \eqref{eq:pressure-energy},
\eqref{eq:pressure-renormalization}, and \eqref{eq:basic-energy} imply
\begin{equation}\label{eq:pressure-damping-lower}
\begin{aligned}
 J^2
 &\le C\left(1+\left(\int\rho^{2\gamma-1}\,dx\right)^2\right)\\
 &\le C\left(1+\int\rho^\gamma\,dx\int\rho^{3\gamma-2}\,dx\right)
 \le C\left(1+\int\frac{aP}{2\mu+\lambda(\rho)}\,dx\right).
\end{aligned}
\end{equation}
H\"older's inequality with \eqref{eq:pressure-renormalization}, \eqref{eq:basic-energy}, and
\eqref{eq:compression-density-moment} implies
\begin{equation}\label{eq:pressure-damping-dual-bound}
\begin{aligned}
 \left\|\frac{a}{2\mu+\lambda(\rho)}\right\|_{L^{2\gamma/(2\gamma-1)}}
 &\le C\|\rho^{2\gamma-2}\|_{L^{2\gamma/(2\gamma-1)}}\\
 &\le C\left(\int\rho^\gamma\,dx\right)^{r-1/(2\gamma)}
       \left(\int\rho^{2\gamma-1}\,dx\right)^{1-r}
 \le C J^{1-r}.
\end{aligned}
\end{equation}
From \eqref{eq:pressure-identities}, \eqref{eq:pressure-damping-lower},
\eqref{eq:pressure-damping-dual-bound}, and \eqref{eq:compression-flux-l2gamma},
\begin{align}
 J'+cJ^2&\le C\bigl(1+A_2(t)+A_3(t)\bigr)J^{1-r},
 \label{eq:pressure-damping}\\
 1+A_2(t)+A_3(t)&\ge cJ^{1+r}\qquad\text{on }\{J'>0\}.
 \label{eq:pressure-growth-threshold}
\end{align}
On $\{J'>0\}$, \eqref{eq:compression-bound} and \eqref{eq:pressure-growth-threshold}
imply
\begin{equation}\label{eq:pressure-positive-variation-bound}
\begin{aligned}
 (J')_+
 &\le C A_3(t)\bigl(1+A_2(t)+A_3(t)\bigr)
       J^{1/2-q_*/(4\gamma)-(1+r)(1-q_*/2)}\\
 &=C A_3(t)\bigl(1+A_2(t)+A_3(t)\bigr)J^{1/2-r}.
\end{aligned}
\end{equation}
On $\{J'\le0\}$, the left-hand side of \eqref{eq:pressure-positive-variation}
vanishes. The exponent equality in \eqref{eq:pressure-positive-variation-bound}
uses \eqref{eq:pressure-variation-exponents} and \eqref{eq:compression-exponent}.
\end{proof}

For $R$ in \eqref{eq:density-maximum}, set
\begin{equation}\label{eq:logarithmic-density-maximum}
 \ell=\log R\ge2,\qquad d\ell=\frac{R'}{R}\,dt.
\end{equation}
\begin{lemma}\label{lem:scalar-cover}
Let $\ell$ be as in \eqref{eq:logarithmic-density-maximum}. For every
$g>0$ and $0<r\le1/2$, we have
\begin{equation}\label{eq:scalar-cover}
 \min\{g,rg^r/\sqrt{\ell}\}+\frac{1}{g\ell^2}\ge\frac{r}{\ell}.
\end{equation}
\end{lemma}
\begin{proof}
For $g\le\ell^{-1}$,
\begin{equation}\label{eq:scalar-cover-small}
 \frac1{g\ell^2}\ge\frac1\ell\ge\frac r\ell.
\end{equation}
For $g\ge\ell^{-1}$,
\begin{equation}\label{eq:scalar-cover-large}
 g\ge\ell^{-1}\ge\frac r\ell,\qquad
 \frac{rg^r}{\sqrt\ell}\ge r\ell^{-r-1/2}\ge\frac r\ell.
\end{equation}
The claim follows from \eqref{eq:scalar-cover-small} and \eqref{eq:scalar-cover-large}.
\end{proof}

\begin{proposition}\label{prop:superquadratic-density}
Let $\gamma>2$ and assume the bounds at $\tau$ in
Proposition~\ref{prop:weighted-flux}. There exists a constant $C_{\tau,T}>0$
such that, for $\tau\le T_1\le T$,
\begin{align}
 \log\log R(T_1)-\log\log R(\tau)
 &\le C_{\tau,T},\label{eq:superquadratic-maximum}\\
 R(T_1)&\le C_{\tau,T}.\label{eq:superquadratic-density-bound}
\end{align}
\end{proposition}
\begin{proof}
With \eqref{eq:pressure-variation-exponents}, set
\begin{equation}\label{eq:maximum-log-variables}
 g=\frac{J}{R^k}.
\end{equation}
Define the time sets
\begin{equation}\label{eq:density-time-sets}
 \mathcal I_J=\{t\in[\tau,T_1]:A_3(t)^2\le J(t)\},\qquad
 \mathcal I_3=\{t\in[\tau,T_1]:A_3(t)^2>J(t)\}.
\end{equation}
From \eqref{eq:energy-comparison},
\begin{align}
 cJ\le\mathcal V\le CJ
 &\quad\text{on }\mathcal I_J,\label{eq:pressure-time-energy}\\
 \mathcal V\le CA_3(t)^2
 &\quad\text{on }\mathcal I_3.\label{eq:dissipation-time-energy}
\end{align}
On $\mathcal I_3$, apply \eqref{eq:logarithmic-density-maximum},
\eqref{eq:maximum-derivatives}, \eqref{eq:dissipation-time-energy},
\eqref{eq:weighted-flux}, and \eqref{eq:basic-energy}.
\begin{equation}\label{eq:dissipation-time-set}
\begin{aligned}
 \int_{\mathcal I_3}\frac{d\ell}{\ell}
 &\le\int_{\mathcal I_3}\frac{\|\widetilde G\|_{L^\infty}}{R\log R}\,dt\\
 &\le C\left(\int_\tau^{T_1}
       \frac{\|\widetilde G\|_{L^\infty}^2}{\mathcal V(R\log R)^2}\,dt\right)^{1/2}
       \left(\int_\tau^{T_1} A_3(t)^2\,dt\right)^{1/2}
 \le C_{\tau,T}.
\end{aligned}
\end{equation}
On $\mathcal I_J$, use \eqref{eq:logarithmic-density-maximum},
\eqref{eq:maximum-log-variables}, \eqref{eq:maximum-derivative-square},
\eqref{eq:pressure-time-energy}, and \eqref{eq:weighted-flux}.
\begin{equation}\label{eq:pressure-time-set}
\begin{aligned}
 \int_{\mathcal I_J}\frac{d\ell}{g\ell^2}
 &=\int_{\mathcal I_J}\frac{R'R^{\gamma-2}}{J\ell^2}\,dt\\
 &\le\int_{\mathcal I_J}\frac{\|\widetilde G\|_{L^\infty}^2}{4J(R\log R)^2}\,dt
 \le C\int_\tau^{T_1}\frac{\|\widetilde G\|_{L^\infty}^2}
                              {\mathcal V(R\log R)^2}\,dt
 \le C_{\tau,T}.
\end{aligned}
\end{equation}
On $\mathcal I_J$, estimate the positive variation of $J$ using
\eqref{eq:pressure-positive-variation}, \eqref{eq:pressure-time-energy},
and Young's inequality.
\begin{equation}\label{eq:variation-pressure-set}
\begin{aligned}
 \frac{rJ^{r-1}(J')_+}{\sqrt{R\log R}}
 &\le C\frac{A_3(t)(1+A_2(t)+A_3(t))}{\sqrt{JR\log R}}\\
 &\le C\left(A_3(t)^2+\frac{(1+A_2(t)+A_3(t))^2}{JR\log R}\right)\\
 &\le C\left(A_3(t)^2+\frac{A_2(t)^2}{\mathcal V R\log R}+1\right).
\end{aligned}
\end{equation}
On $\mathcal I_3$, \eqref{eq:pressure-variation} implies
\begin{equation}\label{eq:variation-dissipation-set}
 \frac{(J')_+}{R^k}\le CA_3(t)\sqrt J\le CA_3(t)^2.
\end{equation}
Define
\begin{equation}\label{eq:variation-primitive}
 \phi(J,R)=\int_0^J\min\!\left\{R^{-k},
                       \frac{rs^{r-1}}{\sqrt{R\log R}}\right\}\,ds .
\end{equation}
From \eqref{eq:variation-primitive}, \eqref{eq:variation-pressure-set},
\eqref{eq:variation-dissipation-set}, \eqref{eq:basic-energy}, and
\eqref{eq:weighted-energies},
\begin{equation}\label{eq:variation-integral}
\begin{aligned}
 \int_\tau^{T_1}\phi_J(J,R)(J')_+\,dt
 &\le\int_{\mathcal I_J}\frac{rJ^{r-1}(J')_+}{\sqrt{R\log R}}\,dt
   +\int_{\mathcal I_3}\frac{(J')_+}{R^k}\,dt\\
 &\le C\int_\tau^{T_1}
   \left(A_3(t)^2+\frac{A_2(t)^2}{\mathcal V R\log R}+1\right)\,dt
 \le C_{\tau,T}.
\end{aligned}
\end{equation}
Substitute $s=R^k y$ into \eqref{eq:variation-primitive} to obtain
\begin{equation}\label{eq:primitive-normalization}
 \phi(J,R)=\Psi(g,\ell),\qquad
 \Psi(g,\ell)=\int_0^g\min\{1,r\ell^{-1/2}y^{r-1}\}\,dy.
\end{equation}
Set
\begin{equation}\label{eq:primitive-threshold}
 g_0=\left(\frac{r}{\sqrt\ell}\right)^{1/(1-r)}.
\end{equation}
Integrate \eqref{eq:primitive-normalization} with threshold
\eqref{eq:primitive-threshold}.
\begin{equation}\label{eq:primitive-branches}
 \Psi(g,\ell)=
 \begin{cases}
  g,&g\le g_0,\\
  \dfrac{g^r}{\sqrt\ell}-\dfrac{1-r}{r}g_0,&g\ge g_0.
 \end{cases}
\end{equation}
The branches in \eqref{eq:primitive-branches} match to first order at $g_0$.
From \eqref{eq:primitive-normalization},
\begin{align}
 \Psi_\ell &\le 0,\label{eq:primitive-ell-derivative}\\
 -R\phi_R &=kg\Psi_g-\Psi_\ell
 \ge k\min\{g,rg^r/\sqrt\ell\}.\label{eq:primitive-r-derivative}
\end{align}
It follows from \eqref{eq:primitive-r-derivative} and
\eqref{eq:variation-integral} that
\begin{equation}\label{eq:primitive-integral}
\begin{aligned}
 \int_\tau^{T_1}\min\{g,rg^r/\sqrt\ell\}\,d\ell
 &\le \frac1k\int_\tau^{T_1}-\phi_R R'\,dt\\
 &=\frac1k\left(\int_\tau^{T_1}\phi_J J'\,dt
      +\phi(J(\tau),R(\tau))-\phi(J(T_1),R(T_1))\right)\\
 &\le \frac1k\left(\int_\tau^{T_1}\phi_J(J')_+\,dt
      +\frac{J(\tau)}{R(\tau)^k}\right)
 \le C_{\tau,T}.
\end{aligned}
\end{equation}
Integrate \eqref{eq:scalar-cover} and use \eqref{eq:pressure-time-set} and
\eqref{eq:primitive-integral} to obtain
\begin{equation}\label{eq:pressure-record-integral}
 r\int_{\mathcal I_J}\frac{d\ell}{\ell}
 \le\int_{\mathcal I_J}\frac{d\ell}{g\ell^2}
    +\int_\tau^{T_1}\min\{g,rg^r/\sqrt\ell\}\,d\ell
 \le C_{\tau,T}.
\end{equation}
\eqref{eq:dissipation-time-set} and \eqref{eq:pressure-record-integral} imply
\begin{equation}\label{eq:superquadratic-maximum-calculation}
 \log\log R(T_1)-\log\log R(\tau)
 =\int_{\mathcal I_J}\frac{d\ell}{\ell}
  +\int_{\mathcal I_3}\frac{d\ell}{\ell}
 \le C_{\tau,T}.
\end{equation}
\eqref{eq:superquadratic-maximum} and the bound on $R(\tau)$ imply
\eqref{eq:superquadratic-density-bound}.
\end{proof}

The two pressure ranges give the finite-time bound used in Section~\ref{sec:4}.

\begin{proposition}\label{prop:finite-density}
For $\gamma>1$, there exists a constant $C_T>0$ such that, for every
$0<T_1\le T$,
\begin{equation}\label{eq:finite-density-bound}
 \sup_{0\le t\le T_1}\|\rho(t)\|_{L^\infty}\le R(T_1)\le C_T.
\end{equation}
The constant depends only on the data in \eqref{eq:basic-regularity}
and is independent of the positive lower bound for $\rho$.
\end{proposition}
\begin{proof}
For $T_1\le\tau_0$, use \eqref{eq:initial-interval}. For $T_1>\tau_0$,
\eqref{eq:density-start-time} and \eqref{eq:initial-interval-start}
allow us to apply \eqref{eq:subquadratic-maximum} for $1<\gamma\le2$
and \eqref{eq:superquadratic-density-bound} for $\gamma>2$.
The dependence of $C_{\tau,T}$ follows from Proposition~\ref{prop:initial-interval}.
\end{proof}

\section{A priori estimates (II): higher-order estimates}\label{sec:4}
\begin{lemma}\label{lem:higher-estimates}
Assume that \eqref{eq:finite-density-bound} holds and that the initial data
satisfy \eqref{eq:initial-hypotheses} for some $q>2$. Then there exists
a constant $C$ depending only on $T$, $q$, $\mu$, $\gamma$,
$\mathbb{T}^2$, $M$, $\|\mathbf{u}_0\|_{H^1}$, and
$\|\rho_0\|_{W^{1,q}}$ such that
\begin{align}
 &\sup_{0\le s\le T_1}
   \bigl(\|\rho(s)\|_{W^{1,q}}+\|\rho_t(s)\|_{L^2}\bigr)
   +\int_0^{T_1}\|\nabla\mathbf{u}\|_{L^\infty}\,ds\le C,
   \label{eq:density-gradient}\\
 &\int_0^{T_1}\!\bigl[\|\mathbf{u}\|_{W^{2,q}}^{1+1/q}
   +s\bigl(\|\mathbf{u}\|_{W^{2,q}}^2
   +\|\nabla\mathbf{u}\|_{L^\infty}^2\bigr)\bigr]\,ds\le C,
   \label{eq:velocity-regularity}\\
 &\sup_{0<s\le T_1}s\|\mathbf{u}(s)\|_{H^2}^2\le C.
   \label{eq:weighted-velocity-h2}
\end{align}
\end{lemma}
\begin{proof}
The proof is similar to that in Huang-Li~\cite{bib3}.
\end{proof}

\section[Proof of the main theorem]{Proof of Theorem~\ref{thm:1.1}}\label{sec:5}
\begin{proposition}\label{prop:5.1}
Assume that \eqref{eq:pressure-law} holds and that $(\rho_0,m_0)$ satisfies
\eqref{eq:local-existence-data}. Then there exists a unique strong solution
$(\rho,\mathbf{u})$ to \eqref{eq:system}--\eqref{eq:initial-data} in
$\mathbb{T}^2\times(0,\infty)$ satisfying \eqref{eq:local-regularity} with $T_0$ replaced by $T$
for every $T>0$. In addition, for any $q>2$, $(\rho,\mathbf{u})$
satisfies \eqref{eq:strong-class} with bounds depending only on $T$, $q$,
$\mu$, $\gamma$, $\mathbb{T}^2$, $M$, $\|\mathbf{u}_0\|_{H^1}$, and
$\|\rho_0\|_{W^{1,q}}$. Moreover, \eqref{eq:basic-regularity} and
\eqref{eq:bounded-density-estimates}--\eqref{eq:bounded-density-weighted-gradient}
hold with a positive constant depending only on $T$, $\mu$, $\gamma$,
$\mathbb{T}^2$, $M$, $\|\mathbf{u}_0\|_{H^1}$, and
$\|\rho_0\|_{L^\infty}$.
\end{proposition}
\begin{proof}
The proof is similar to that of Huang-Li~\cite{bib3}
and is omitted.
\end{proof}
\noindent\textbf{Proof of Theorem~\ref{thm:1.1}.}
Let $(\rho_0,m_0)$ satisfying \eqref{eq:initial-hypotheses} be the initial
data as described in Theorem~\ref{thm:1.1}. For constant $\delta\in(0,1)$,
we define
\begin{equation}\label{eq:positive-initial-approximation}
 \rho_0^\delta=\frac{M}{M+\delta}(j_\delta*\rho_0+\delta)
 \ge \frac{M\delta}{M+\delta}>0,\qquad
 \mathbf{u}_0^\delta=j_\delta*\mathbf{u}_0,\qquad
 m_0^\delta=\rho_0^\delta\mathbf{u}_0^\delta,
\end{equation}
where $j_\delta$ is the standard periodic mollifying kernel of width
$\delta$. Hence, we have $\rho_0^\delta,\mathbf{u}_0^\delta\in H^\infty$, and
\begin{equation}\label{eq:initial-approximation-convergence}
 \lim_{\delta\to0}\left(
 \|\rho_0^\delta-\rho_0\|_{W^{1,q}}
 +\|\mathbf{u}_0^\delta-\mathbf{u}_0\|_{H^1}\right)=0.
\end{equation}

Proposition~\ref{prop:5.1} thus yields that the problem
\eqref{eq:system}--\eqref{eq:initial-data} with $(\rho_0,m_0)$ being
replaced by $(\rho_0^\delta,m_0^\delta)$ has a unique global strong
solution $(\rho^\delta,\mathbf{u}^\delta)$ satisfying
\eqref{eq:strong-class} for any $T>0$, with bounds independent of $\delta$.
Moreover, there exists a positive constant $C_T$ depending only on $T$,
$\mu$, $\gamma$, $\mathbb{T}^2$, $M$, $\|\rho_0\|_{L^\infty}$, and
$\|\mathbf{u}_0\|_{H^1}$ such that $(\rho^\delta,\mathbf{u}^\delta)$
satisfies \eqref{eq:basic-regularity}. Letting $\delta\to0$, standard
arguments (see \cite{bib3}) thus show that
\eqref{eq:system}--\eqref{eq:initial-data} has a global strong solution
$(\rho,\mathbf{u})$ satisfying \eqref{eq:strong-class}. Moreover,
$(\rho,\mathbf{u})$ satisfies \eqref{eq:basic-regularity} with the same
constant $C_T$.

Since the proof of the uniqueness of $(\rho,\mathbf{u})$ satisfying
\eqref{eq:strong-class} is similar to that of Germain~\cite{Germain2011}.
\section*{Acknowledgments}
	X. Huang is partially supported by Chinese Academy of Sciences Project for Young Scientists in Basic Research (Grant No. YSBR-031), National Natural Science Foundation of China (Grant Nos. 12494542, 11688101) and National Key R\&D Program of China (Grant No. 2021YFA1000801). 
	\vspace{0.3cm}
    
\subsection*{AI assistance}

The authors acknowledge the use of AI to explore proof strategies,
check calculations, and polish the manuscript. In particular, AI
suggested the double-commutator argument used to obtain additional
density integrability at $\beta=1$ and helped verify the cancellation
in the pressure correction for $\gamma>2$. The authors take full
responsibility for the mathematical arguments and the final manuscript.
\vspace{0.3cm}

	\noindent\textbf{Data availability statement.} Data sharing is not applicable to this article.
	
	\vspace{0.3cm}
	\noindent\textbf{Conflict of interest.} The authors declare that they have no conflict of interest.
\phantomsection

\end{document}